\documentclass[11pt]{article}
\def\be{\begin{eqnarray}}
\def\ee{\end{eqnarray}}

\usepackage{graphics}
\usepackage{subfigure}
\usepackage{amsmath}
\usepackage{amsfonts}
\usepackage{amsbsy}
\usepackage{epsfig}
\usepackage{harvard}

\title{\textbf{Non local Lotka-Volterra  cross-diffusion systems in heterogeneous media}}

\def\me{\medskip\noindent}
\def\bi{\bigskip\noindent}
\def\be{\begin{eqnarray}}
\def\ee{\end{eqnarray}}
\def\ben{\begin{eqnarray*}}
\def\een{\end{eqnarray*}}

\newtheorem{theo}{\bf Theorem}[section]
\newtheorem{lem}[theo]{\bf Lemma}
\newtheorem{prop}[theo]{\bf Proposition}
\newtheorem{coro}[theo]{\bf Corollary}
\newtheorem{rem}[theo]{\bf Remark}
\def\RR{\mathbb{R}}
\def\NN{\mathbb{N}}
\def \EE{\mathbb{E}}
\def \PP{\mathbb{P}}

\def \1{\mathbf{1}}

\renewcommand{\Box}{\hfill\mbox{\fbox{\rule{0mm}{1.5mm}}}}

\author{Joaquin Fontbona\thanks{
DIM-CMM, UMI(2807) UCHILE-CNRS, Universidad de Chile, Casilla 170-3, Correo 3,
Santiago-Chile, E-mail: \texttt{fontbona@dim.uchile.cl}.  } \,  and  Sylvie
  M\'el\'eard\thanks{CMAP, Ecole Polytechnique, CNRS UMR 7641, route de Saclay, 91128 Palaiseau
    Cedex-France; E-mail: \texttt{sylvie.meleard@polytechnique.edu}. }}

\date{\today}

\begin{document}

\maketitle

\pagenumbering{roman} \pagestyle{plain} \thispagestyle{empty}

\pagenumbering{arabic}

\begin{abstract}

We introduce a  stochastic individual model  for the spatial behavior of  an animal  population of dispersive and competitive species, considering various kinds of biological effects, such as heterogeneity of environmental conditions, mutual attractive or repulsive interactions between individuals or  competition between them  for resources.  As a consequence of the  study of the large population limit,   global existence of a nonnegative weak solution to a multidimensional  parabolic strongly coupled model of competing species is proved. The main new feature of the corresponding integro-differential equation is the nonlocal nonlinearity appearing in the diffusion terms, which may  depend on the spatial densities of all population  types.  Moreover, the diffusion matrix is generally not strictly positive definite and  the cross-diffusion effect allows for   linearly  growing influences  of the subpopulations' sizes.  We prove  uniqueness of the finite measure-valued solution and give conditions under which the solution takes values in a functional space. We then make the competition kernels converge to a Dirac measure and obtain the existence of a solution to  a locally competitive version of the previous equation. The large population limit is obtained by means of  weak convergence tools for  measure valued processes.  The techniques  employed in the study of the liming equation are novel, and are essentially based on the underlying stochastic flow related to the dispersive part of the dynamics, together with the use of suitable dual distances in the  space of finite measures. 
\smallskip 

This is an extended version of the paper ``Non local Lotka-Volterra system with cross-diffusion in a heterogeneous medium'' appeared in J. Math. Biol (2015) 70:829-854.
\end{abstract}

\section{Introduction}

 The spatial structure of a biological community and  in particular the spatial distribution formed by dispersive motions of populations with intra- and inter- specific interactions  is a fundamental subject in  mathematical ecology \citeaffixed{KS70,GN76,GN75,MM,MY,SKT}{see e.g.}. 
In the present work, the  spatial behavior of  a population of competitive species is studied. The dispersive motion of an individual in its environment is modeled as the result of various kinds of biological effects, such as heterogeneity of environmental conditions, mutual attractive or repulsive interactions with  other  individuals and  competition  for resources.  These different effects will be modeled by local interaction kernels depending on the type of the individual and acting either on its spatial parameters or on its ecological parameters. 

\me
The population is composed of $M$ sub-populations (species) characterized by different phenotypes. Each species has its own spatial and ecological dynamics depending on the spatial and genetic characteristics of the whole population. We assume that the motion of each individual (of a given type) is driven by a diffusion process on $\mathbb{R}^d$ whose coefficients depend on the spatial repartition of the different species around. Moreover, the individuals may reproduce and die, either from their natural death or because of the competition pressure for sharing  resources. Each species has its own growth rate. The competition pressure of an individual of type $j$ on an individual of type $i$ depends both on the location of these  individuals and on their type. It is called intra-specific competition if $i=j$,  and inter-specific in case $i\neq j$.  It is not assumed to be symmetric in $(i,j)$.

We  describe the stochastic dynamics of such a population by an individual-based model. Each individual is characterized by its type and its spatial location. Because of the births and deaths of individuals, the population doesn't live in a vector space of positions. We model its dynamics  as a Markov process  with   values in the $M$-dimensional vector space of   $\mathbb{R}^d$-point measures.

We introduce the charge capacity parameter $K$ describing  the order of the population size so that, to  be consistent, the individuals are weighted by ${1\over K}$. The existence of the population process is obtained by standard arguments. Then,  large population asymptotics  is studied, using classic tools for the convergence of measure valued processes.  

We  show that when $K$ tends to infinity, the population process converges to a weak solution of the following nonlocal (in trait and space)  nonlinear parabolic cross-diffusion-reaction  system: for all $i\in \{1, \cdots, M\}$,
\be 
\partial_{t}u^i &=&\frac{1}{2} \sum_{k,l=1}^d \partial^2_{x_{k}x_{l}} \left(a^i_{k,l}(.,G^{i1}*u^1, \cdots, G^{iM}*u^M) \,u^i\right)\nonumber\\
&& - \sum_{k=1}^d\partial_{x_{k}} \left(b^i_{k}(.,H^{i1}*u^1, \cdots, H^{iM}*u^M) \,u^i\right) + \left( r_{i} - \sum_{j=1}^M C^{ij}*u^j\right) \, u^i.
\label{cross-diffusion1}
\ee
Here,   $u^i(t,.)$ denotes  general a finite measure on $\mathbb{R}^d$ for any $t\geq 0$, and   $(G^{ij}, H^{ij}, C^{ij})_{1\leq i,j\leq M}$  are  $3M^2$ nonnegative and smooth functions defined from $\mathbb{R}^d$ to $\mathbb{R}_{+}$  that  model the spatial  interactions between  individuals of type $i$ and $j$. 

By  means of this convergence result,  we get a theorem of  existence of a weak solution to Equation \eqref{cross-diffusion1}.
Next we prove the uniqueness of such a solution and we give two sets of assumptions under which the measure solution has a density with respect to the Lebesgue measure, thus establishing the existence of a function solution to \eqref{cross-diffusion1}. The tools we use to establish these  and the forthcoming results are again probabilistic ones. 

\me
In the model leading to Equation  \eqref{cross-diffusion1},  the competition between two individuals  is described as a function of the distance between them.  This biological assumption is clear: the closer  the animals are, the stronger is the fight to share resources. An extreme situation is the local case, when individuals only compete if they stay at the same place.  Mathematically speaking, this means   that for $i,j\in \{1, \cdots, M\}$ and for $x\in \mathbb{R}^d$, the competition kernel has the form  
$$C^{ij}(x) = c_{ij} C_{\varepsilon}(x),$$ 
where $c_{ij}$ are positive constant numbers and the measures  $C_{\epsilon}(x-y) dy$   weakly converge to the Dirac measure at $x$ when the range of  interaction $\varepsilon$ tends to $0$ (for instance,   $ C_{\varepsilon}$ may be the centered Gaussian density with variance $\varepsilon$).  In this setting, 
we  study the convergence of the solution $u^\varepsilon$  of \eqref{cross-diffusion1} when $\varepsilon$ tends to zero, showing its  convergence to the unique solution $u$ of the spatially nonlocal nonlinear cross-diffusion equation: for all $i\in \{1, \cdots, M\}$,
\be 
\partial_{t}u^i &=& \frac{1}{2}\sum_{k,l=1}^d \partial^2_{x_{k}x_{l}} \left(a^i_{k,l}(.,G^{i1}*u^1, \cdots, G^{iM}*u^M) \,u^i\right)\nonumber\\
&&  - \sum_{k=1}^d\partial_{x_{k}} \left(b^i_{k}(.,H^{i1}*u^1, \cdots, H^{iM}*u^M) \,u^i\right) + \left( r_{i} - \sum_{j=1}^M c_{ij}\, u^j\right)  u^i.
\label{cross-diffusion2} 
\ee
 Ecological models featuring 
space displacements have been studied by \citeasnoun{CM07},  \citeasnoun{ADP} and  \citeasnoun{B},  but  the diffusion coefficients therein only depend on the type of each  individual, and not on the spatial distribution of the other animals alive.  
To our knowledge, the nonlocal nonlinear Equations \eqref{cross-diffusion1}  and  \eqref{cross-diffusion2} have never been studied in such generality, despite the fact that they  naturally arise from the biological motivation. A recent paper on  conservative relaxed cross-diffusion  \citeasnoun{LPR} addresses  well posedness issues in a    model with nonlocal interaction in the diffusion terms. Nonlocal reaction terms have  otherwise been considered  by  \citeasnoun{CD},   \citeasnoun{BNPR},  \citeasnoun{GVA} and in references therein. Our model incorporates those two features simultaneously.

The dependance of the  diffusion coefficient on  the individual density is indeed the main new difficulty we deal with. It is the reason why the introduction of different techniques from those developed in the aforementioned works will be needed.  

Cross-diffusion models with local spatial interaction and local competition have excited the scientific community, see for example the works of  \citeasnoun{MM},  \citeasnoun{MK},  \citeasnoun{LMN00},  \citename{CJ04} \citeyear{CJ04,CJ06} and more recently in a work by \citeasnoun{DLM} by developing suitable entropy methods.  The prototypical equation in this situation has the form
 \begin{equation}\label{functionsol}
\begin{split}
\partial_{t} u^1 =  &d_{1} \Delta \left( (a_{1} + b_{12}u^2)u^1 \right) + (r_{1} - c_{11} u^1 - c_{12} u^2)u^1\, , \\
 \partial_{t} u^2 = &d_{2} \Delta \left( (a_{2} + b_{21}u^2)u^1 \right)   + (r_{2} - c_{21}u^1 - c_{22} u^2)u^2\, , \\ 
\end{split}
 \end{equation}
 with $a_i\geq 0$, $b_{ij}>0$, $i,j=1,2$  and boundary conditions on a given bounded smooth domain of $\mathbb{R}^d$.  Global existence results were obtained  by  \citename{CJ04} \citeyear{CJ04,CJ06}  for such  equation with $a_i>0$. We refer the  works \cite{DLMT}  and \cite{Jung2015} for more recent  (and generalized) developments based on entropy approaches.  We notice that those techniques do  not apply, in principle, to  nonlocal interactions at the level of the diffusion and drift coefficients as we will consider here. In turn, using our techniques, we cannot  yet recover well possedned results for equations like \eqref{functionsol}. Establishing the convergence of the system  \eqref{cross-diffusion2}  to systems of the type of  
 \eqref{functionsol} when the kernels $G$ and $H$ tend to Dirac measures,  is indeed to be a highly difficult  open  challenge.

\medskip

For later use, let us  introduce some notation:

\begin{itemize}
\item We denote by ${\cal M}$ the space of finite measures in $\mathbb{R}^d$ endowed with the weak topology. 
\item  For $k\geq 1$ and $\alpha \in (0,1)$, we denote by ${\cal C}^{k,\alpha}(\RR^d)$ the space of $k$ times  differentiable  functions  on $\RR^d$  which have  bounded derivatives up to the  $k-$th   order and a  globally  $\alpha$-Holder derivative of order $k$. Notice that   functions  in ${\cal C}^{k,\alpha}(\RR^d)$ are not required to be bounded. 
\item  The subspace of bounded functions in ${\cal C}^{k,\alpha}(\RR^d)$ is denoted by ${\cal C}^{k,\alpha}_b(\RR^d)$. 
\item  The notation  ${\cal C}^{k}(\RR^d)$ and  ${\cal C}^{k}_b(\RR^d)$ is defined analogously, without the Holder continuity requirement. 
\item  The super-indexes in the above notations will be omitted  when $k=\alpha=0$. 
\item ${\cal C}_0(\RR^d)$ denotes the space of compactly supported continuous functions on $\RR^d$. 
\item ${\cal C}^{1,2}([0,T]\times \mathbb{R}^d)$ will denote the space of functions on $[0,T]\times \mathbb{R}^d$ which are continuously differentiable up to the order $2$ in the space variable and continuously differentiable in the time variables  (with lateral limits when needed).
\end{itemize}

\section{The Individual-based Model}
\subsection{Assumptions}\label{Assumptions}
 Let us denote by ${\cal S}_+(\mathbb{R}^d)$  the space of symmetric nonnegative diffusion matrices and   define for $i=1,\cdots, M$ the measurable functions
\ben && a^i: \mathbb{R}^d\times \mathbb{R}_+^M \to {\cal S}_+(\mathbb{R}^d), \\
&& b^i: \mathbb{R}^d\times \mathbb{R}_+^M \to \mathbb{R}^d , \\
&& r_{i} : \mathbb{R}^d \to \mathbb{R}_{+}. \een
We denote by $\sigma^i$ the $d\times d-$matrix such that $a^i=\sigma^i\,(\sigma^i)^*$. We will assume throughout this work the following hypotheses:
\medskip

$(\mathbf{H})$:
\begin{itemize}

\item[i)] There is a positive constant $L$ such that for any $x,x'\in \RR^d$ and  any $ v_{j},v'_{j}\in \RR_+ $  with $ j\in \{1,\cdots, M\},$
\ben
&&|\sigma^i(x,v_{1},\cdots, v_{M})-\sigma^i(x',v'_{1},\cdots, v'_{M})|+|b^i(x,v_{1},\cdots, v_{M})-b^i(x',v'_{1},\cdots, v'_{M})|\\
&&\leq L \Big(|x-x'|+\sum_{j=1}^M|v_{j}-v'_{j}|\Big).\een
Moreover, there exists $C_{\sigma}>0$ such that for all $x\in \RR^d$ and $v=(v_{1},\cdots, v_{M})\in \RR_+^M$,
 $$ |\sigma^i(x,v_{1},\cdots, v_{M})|\leq C_{\sigma} (1+ |v|).$$
\item[ii)] The  functions  $(G^{ij}, H^{ij}, C^{ij})_{1\leq i,j\leq M}$  defined from $\mathbb{R}^d$ to $\mathbb{R}_{+}$ are assumed  to be nonnegative, bounded and Lipschitz continuous.
\item[iii)] The nonnegative functions $r_{i}$ are assumed to be  continuos and bounded by $\bar{r}>0$.
\end{itemize}

%

\subsection{The  diffusive $M$-type  stochastic population dynamics}

Let us now describe the dynamics of the population we are interested in.    The population dynamics  will be modeled by a point measure-valued Markov processes, which takes into account the births and deaths of all individuals, as well as their motion during their life.

 Let us   fix the charge capacity $K\in \mathbb{N}^*$ and define
$${\cal M}_K =\left\{ {1\over K} \sum_{n=1}^N \delta_{x^n}, x^n\in \mathbb{R}^d,N\in \NN\right\}$$
 as the space of weighted finite point measures on $\mathbb{R}^d$.
  The  stochastic population process $(\nu^K_t)_{t\geq 0}$ representing the populations' dynamics will
take values in $({\cal M}_{K})^M$. The $i$th coordinate of this process describes the spatial configuration of the subpopulation of type $i$. Thus,
$$\nu^K_t = ( \nu^{1,K}_t, \cdots,  \nu^{M,K}_t)=( \nu^{i,K}_t)_{1\leq i\leq M}$$
with
\begin{eqnarray}\label{processnu}
  \nu^{1,K}_t = {1\over K}\sum_{n=1}^{N^1_t} \delta_{X^{n,1}_t}\  \  ,\cdots\ ,\  \nu^{M,K}_t ={1\over K}\sum_{n=1}^{N^M_t} \delta_{X^{n,M}_t},
\end{eqnarray}
where for any $i\in \{1,\cdots,M\}$, $N^i_t =K \langle \nu_i^{i,K},1 \rangle \in \NN$  stands for the number of living individuals of type $i$ at time
$t$ and $X^{1,i}_t,\dots ,X^{N^i_t,i}_t$ indicate  their positions in $\mathbb{R}^d$.
  
   \bi

The dynamics of the population can be roughly summarized as follows:

\me
\begin{itemize}
\item  The initial population is characterized by the measures $(\nu_0^{i,K})_{1\leq i\leq M} \in ({\cal
  M}_{K})^M$ at time $t=0$. Any individual of  type $i$ located at $x \in \mathbb{R}^d$ at time $t$ has two independent exponential
clocks: a ``clonal reproduction''  clock with parameter $r_{i}(x)$  and a
``mortality'' clock with parameter $\sum_{j=1}^M{C^{ij}}*  \nu^{j,K}_t(x)  $. If the
reproduction  clock of an individual rings, then it produces
at the same location an individual  of same type as itself. If its mortality
clock  rings, then the individual disappears. The  death rate of an $i$th  type individual depends on the positions of  the other individuals through the kernel $C^{ij}$, which describes how species $j$ acts on $i$ in the competition for resources.

\item 
During its life, an individual will move as a diffusion process whose coefficients depend on all individual positions.
The motion of an individual with type $i$ is a  diffusion process with diffusion
matrix $a^i(.,G^{i1}*\nu_t^{1,K}, \cdots, G^{iM}*\nu_t^{M,K}) $  and drift vector $b^i(.,H^{i1}*\nu_t^{1,K}, \cdots, H^{iM}*\nu_t^{M,K})$.  The coefficients thus  take into account the  non-homogeneous spatial densities of the  different species:  a species can be attracted or repulsed by the other ones and the concentration of species  may increase or decrease the fluctuations in the dynamics. 

\end{itemize}


The vector-measure valued process $(\nu^K_t)_{t\geq 0}$ will be a Markov process that can be rigorously constructed as solution of a stochastic differential equation driven by  $d$-dimensional Brownian motions $(B^{n,i})_{{1\leq i \leq M, n \in \mathbb{N}^*}}$ and  Poisson point measures $(Q^i(dt, dn, d\theta))_{1\leq i \leq M}$ on $\mathbb{R}_{+}\times \mathbb{N}^* \times \mathbb{R}_{+}$ with intensity $dt\otimes \sum_{m\in \mathbb{N}^*} \delta_{m}(dn)\otimes d\theta$, all independent and independent of the initial condition $(\nu^{1,K}_{0},\dots, \nu^{M,K}_{0})$. 

Let us give an explicit algorithmic construction of  $(\nu^K_t)_{t\geq 0}$  in terms of these objects.   To that end, we fix first an arbitrary ordering $\preceq$ in $\RR^d$ (say, the lexicographical one), which will be used to define a numbering of the positions of the individuals  $(X^{1,i}_t,\dots ,X^{N^i_t,i}_t)$ of  the population $i$ at some specific time instants $t\geq 0$.  We then proceed as follows: 

\begin{itemize} 
\item[0.]  Set $m=0$  and $T_0=0$.  For each $i\in\{1,\dots ,M\}$, the vector  $(X^{1,i}_0,\dots ,X^{N^i_0,i}_0)$ is defined as the $\preceq$ -ordered elements of the support of $\nu^{i,K}_0$.  
\item[1.] For each $i\in\{1,\dots ,M\}$ the $i$th  population follows the dynamics 
\begin{equation}\label{ithparticle}
\begin{split}
X^{n,i}_t=X^{n,i}_{T_m}+ &\int_{T_m}^t \sigma^i(X^{n,i}_{s}, G^{i1}*\nu^{1,K}_{s}(X^{n,i}_{s}), \cdots, G^{iM}*\nu^{M,K}_{s}(X^{n,i}_{s})) dB^{n,i}_s\\
 +&\int_{T_m}^t b^i(X^{n,i}_{s}, H^{i1}*\nu^{1,K}_{s}(X^{n,i}_{s}), \cdots, H^{iM}*\nu^{M,K}_{s}(X^{n,i}_{s})) ds, \, \, n=1,\dots, N^i_{T_m} \\
\end{split}
\end{equation}
until  the next jump time $t$ of any of the measures  
$$Q^j(dt,dn,d\theta), \quad  j\in \{1,\dots, M\},$$
such that the  corresponding atom $(t,n,\theta)\in \mathbb{R}_{+}\times \mathbb{N}^* \times \mathbb{R}_{+}$ and index $j\in\{1,\dots ,M\}$ satisfy 
\begin{equation}\label{jumps}
 n\leq  N^j_{T_m}\,    \mbox{  and   }\, \theta\leq r_{j}(X^{n,j}_{t-})+\sum_{l=1}^M C^{jl} *\nu^{l,K}(X^{n,j}_{t-}) .
\end{equation}
We then set $T_{m+1}= t$. 
\item[2.]  For $i=j$ and $(t,n,\theta)=(T_{m+1},n,\theta)$ as before, 
\begin{itemize}
\item if  $\theta\leq r_{i}(X^{n,i}_{t-})$ we set $N^i_{t}=N^i_{t-}+1$ and define the vector  $$\left(X^{1,i}_{t},\dots ,X^{N^i_{t},i}_{t}\right)$$
as the increasing  $\preceq$ -rearrangement of $(X^{1,i}_{t-},\dots ,X^{N^i_{t-},i}_{t-}, X^{n,i}_{t-})$ or, otherwise,
\item if  $ r_{i}(X^{n,i}_{t-})< \theta\leq r_{i}(X^{n,j}_{t-})+\sum_{l=1}^M C^{il} *\nu^{l,K}(X^{n,i}_{t-})  $ we set $N^i_{t}=N^i_{t-}-1$ and define the vector  $$\left(X^{1,i}_{t},\dots ,X^{N^i_{t},i}_{t}\right)$$
as the increasing  $\preceq$ -rearrangement of  $(X^{1,j}_{t-},\dots , X^{n-1,j}_{t-},  X^{n+1,j}_{t-}\dots, X^{N^j_{t-},j}_{t-})$.
\end{itemize} 

For $i\neq j$, we take $\left(X^{1,j}_{t},\dots ,X^{N^j_{t},j}_{t}\right)= \left(X^{1,j}_{t-},\dots ,X^{N^j_{t-},j}_{t-}\right)$. 

\item[3.] We increase $m$ by $1$ and go to  step $1.$
\end{itemize}


We have
\begin{lem} Assume that  $\EE(\langle \nu_0^{i,K}, 1\rangle) <\infty$ for $i=1,\dots. M$.  Then, 
the process $ ( \nu^{i,K}_t)_{1\leq i\leq M} $  given by  \eqref{processnu} is defined for all $t\geq 0$.
\end{lem}

{\bf Proof:}   Notice first  that $(T_m)_{m\in \NN}$ are stopping times with respect to the filtration $({\cal F}_t)_{t\geq 0}$ generated by  the processes $(B^{n,i})_{{1\leq i \leq M, n \in \mathbb{N}^*}}$ and   $(Q^i(dt, dn, d\theta))_{1\leq i \leq M}$ on $\mathbb{R}_{+}\times \mathbb{N}^* \times \mathbb{R}_{+}$  and the measures  $(\nu_0^{1,K},\dots, \nu^{M,K}_{0})$. Moreover, thanks to $(\mathbf{H})$,   conditionally  on ${\cal F}_{T_m}$  the  coefficients of  equations \eqref{ithparticle} seen as an  SDE in $\prod_{i=1}^M \RR^{ d \, N_{T_m}^{i}}$, are globally Lipschitz continuous (with Lipschitz constants depending on $(N_{T_m}^{1},\cdots,N_{T_m}^{M})$) and so the processes  $(X_t^{n,i})_{1\leq n\leq  N^i_ {T_m}, 1\leq i \leq M}$ can be defined,  for  all $t\in [T_m,T_{m+1})$,   as the restriction to that interval of some process of fixed dimension defined on $[T_m,\infty)$.

The process $ ( \nu^{i,K})_{1\leq i\leq M} $  is thus  defined  on $[0, T_{\infty})$ with  $T_{\infty}:=\lim_{m\to \infty} T_m $ so we need to prove that $T_{\infty} =\infty$ a.s.  Writing $S_R:=\inf\{s\in [0, T_{\infty}) : \langle \nu_s^{i,K}, 1\rangle \geq R\}$ for each $R>0$, by construction  of $\nu^{i,K}$ we get
\begin{equation}\label{mass}
\begin{split}
& \langle \nu^{i,K}_{t\wedge S_R}, 1\rangle=  \langle \nu^{i,K}_{0}, 1\rangle \\
& + 
\int_{[0,t\wedge S_R]\times \mathbb{N}\times  \mathbb{R}_{+}}  {\left( {\bf 1}_{ \theta\leq r_{i}(X^{n,i}_{s-})} -  {\bf 1}_{r_{i}(X^{n,i}_{s-})< \theta\leq r_{i}(X^{n,i}_{s-})+\sum_{j=1}^M C^{ij} *\nu^{j,K}(X^{n,i}_{s-})}\right)   \over K}   {\bf 1}_{ n\leq  K\,\langle \nu^{i,K}_{s-},1\rangle}  Q^i(ds,dn,d\theta),
\\
\end{split}
\end{equation}
hence
$$\sup_{t\in [0,T\wedge S_R]} \langle \nu^{i,K}_{t}, 1\rangle\leq   \langle \nu^{i,K}_{0}, 1\rangle
 + 
{ 1\over K} \int_{[0,T\wedge S_R]\times \mathbb{N}\times  \mathbb{R}_{+}}     {\bf 1}_{  \theta\leq \bar{r}, n\leq  K\,\langle \nu^{i,K}_{s-},1\rangle}  Q^i(ds,dn,d\theta)$$
and, by compensation, 
\begin{equation}\label{bycompen}
\begin{split}
 \EE\left(\sup_{t\in [0,T\wedge S_R]} \langle \nu^{i,K}_{t}, 1\rangle\right)\leq  & \,  \EE\left(  \langle \nu^{i,K}_{0},1\rangle\right)
 + \bar{r}\int_0^T \,  \EE\left(\langle \nu^{i,K}_{s\wedge S_R}, 1\rangle \right)ds	\\
 \leq & \,  \EE\left( \langle \nu^{i,K}_{0},1\rangle \right)
 + \bar{r}\int_0^T \EE\left( \sup_{s\in [0,t\wedge S_R]} \langle \nu^{i,K}_{s}, 1\rangle \right)dt.	\\
\end{split}
\end{equation}
Gronwall's lemma  then yields $ \EE\left(\sup_{t\in [0,T\wedge S_R]} \langle \nu^{i,K}_{t}, 1\rangle\right)\leq e^{\bar{r} T}  \EE\left(  \langle \nu^{i,K}_{0},1\rangle\right)
$ and, by monotone convergence as $R\nearrow \infty$,     $ \EE\left(\sup_{t\in [0,T\wedge T_{\infty}]} \langle \nu^{i,K}_{t}, 1\rangle\right)<\infty$ for all $i\in \{1,\dots,M\}$. Since by  \eqref{jumps} the number of jumps of  $\nu^{i,K} $  on the interval $ [0,T\wedge T_{\infty}]$ is bounded by 
$$\int_{[0,T\wedge T_{\infty}]\times \mathbb{N}\times  \mathbb{R}_{+}}     {\bf 1}_{  \theta
 \leq \bar{r}+ K \sum_{l=1}^M \| C^{il} \|_{\infty}\sup_{t\in [0,T\wedge T_{\infty}]} \langle \nu^{l,K}_{t}, 1\rangle   ,\,  n\leq  K\,\sup_{t\in [0,T\wedge T_{\infty}]} \langle \nu^{i,K}_{t}, 1\rangle}  Q^i(ds,dn,d\theta)<\infty,$$ 
we deduce that $T_{\infty}\geq T$ for each $T>0$ and this completes the proof.

\Box

\begin{rem}
It is clear from their construction that the real processes  $t\mapsto \langle \nu^{j,K}_t ,f \rangle= {1\over K}\sum_{n=1}^{N^j_t} f(X^{n,j}_t)  $    are right continuous with left limits for each $f\in {\cal C}_b(\RR^d)$, hence  the processes $t\mapsto   {1\over K}\sum_{n=1}^{N^j_t} \delta_{X^{n,j}_t}$  are right continuous with left limits in the space ${\cal M}$ endowed with the weak topology.  
\end{rem}

For each $\nu=(\nu^i)_{1\leq i\leq M} \in {\cal M}^M$,  let us now  introduce  a  second order differential operator acting on real  $ {\cal C}^2(\RR^2)$   functions $f$   by
\begin{equation*}
\begin{split}  {\cal L}^i_{\nu} f(x):= &  Tr  (\big(  a^i (x,G^{i1}*\nu^{1}(x), \cdots, G^{iM}*\nu^{M}(x)) \mbox{Hess} f(x)) \big) \\
& + b^i(x,H^{i1}*\nu^{1}(x), \cdots, H^{iM}*\nu^{M}(x))\cdot \nabla f(x)\\
\end{split}
\end{equation*}

By  applying It\^o's formula to $\left(X^{1,i}_{s},\dots ,X^{N^i_{T_m},i}_{s}\right)$ for $s\in [T_m,T_{m+1})$, summing over $n\in \{1,\dots, N^i_{T_m}=K \langle \nu^i_{s},1\rangle\}$ and then, for  given $t\geq 0$,  over $m\in \NN$  such that $T_m\leq t$, the process $(\nu_t ^{i,K})_{t\geq 0} $  is seen to satisfy   for  any  function $f\in{\cal C}^{1,2}(\RR_+\times \mathbb{R}^d)$ the stochastic differential equation: 
\begin{equation}\label{ito}
\begin{split}
\langle  \nu& ^{i,K}_{t},  f(t,\cdot) \rangle=  \langle \nu^{i,K}_{0}, f(0,\cdot)\rangle  + 
 \int_{0}^t 
\bigg\langle \nu^{i,K}_{s}, \,  \partial_s f (s,\cdot)+ {\cal L}^i_{\nu_s^K}f(s,\cdot) \bigg\rangle ds \\
 &  +  {1\over K}\int_{0}^t \sum_{n=1}^{K \langle \nu^i_{s},1\rangle} \nabla^* f(s,X^{n,i}_{s})   \sigma^i(X^{n,i}_{s}, G^{i1}*\nu^{1,K}_{s}(X^{n,i}_{s}), \cdots, G^{iM}*\nu^{M,K}_{s}(X^{n,i}_{s})) ) \,dB^{n,i}_s \\
& + 
\int_{[0,t]\times \mathbb{N}\times  \mathbb{R}_{+}} \left( {\bf 1}_{\theta\leq r_{i}(X^{n,i}_{s-})} -  {\bf 1}_{r_{i}(X^{n,i}_{s-})< \theta\leq r_{i}(X^{n,i}_{s-})+\sum_{j=1}^M C^{ij} *\nu^{j,K}(X^{n,i}_{s-})}\right)  \\  &\hspace{6cm}\times { f(s-,X^{n,i}_{s-}) \over K} {\bf 1}_{ n\leq  K\,\langle \nu^i_{s-},1\rangle}  Q^i(ds,dn,d\theta).\\
\end{split}
\end{equation}

By compensation, we can rewrite
\begin{equation}\label{itocomp}
\begin{split}
\langle  \nu^{i,K}_{t},  f(t,\cdot)\rangle=   &\langle \nu^i_{0}, f(0,\cdot)\rangle  + 
 \int_{0}^t 
\bigg\langle \nu^{i,K}_{s}, \partial_s f (s,\cdot)+ {\cal L}^i_{\nu_s^K}f (s,\cdot) + \big(r_{i}- \sum_{j=1}^M C^{ij}* \nu^{j,K}_{s}\big)f (s,\cdot) \bigg\rangle ds \\
& + M^{i,f,B}_t+ M^{i,f,Q}_t \,  , \\
\end{split}
\end{equation}
where  $M^{i,f,B}$ and $ M^{i,f,Q}$ are, respectively, a continuous  local martingale with quadratic variation given by 
\begin{equation}\label{quadvarB}
\langle M^{i,f,B} \rangle_t= \frac{1}{K} \int_0 ^t \big\langle\nu^{i,K}_s,    \nabla^* f (s,\cdot) a^i (\cdot,G^{i1}* \nu^{i,K},\cdots,G^{iM}* \nu^{M,K})  \nabla f(s,\cdot)   \big\rangle ds
\end{equation}
and a compensated  pure jump  local martingale  with predictable  quadratic variation
\begin{equation}\label{quadvarQ}
  \langle  M^{i,f,Q} \rangle_t= \frac{1}{K}\int_{0}^t 
\big\langle \nu^{i,K}_{s}, \big(r_{i}+ \sum_{j=1}^M C^{ij}* \nu^{j,K}_{s}\big)f^2 (s,\cdot) \big\rangle ds . \end{equation}
The local martingales $M^{i,f,B}$ and $ M^{i,f,Q}$ will be fundamental in the study of the limit of   $(\nu^K_{t})_{t\in [0,T]}$ when $K$ goes to $\infty$. 
%

\medskip 

Using the above formulae we can  characterize the law of  $(\nu^K_{t})_{t\in [0,T]}$  as a process with values in ${\cal M}^M$.  To that end, observe first  that the zero order term on the r.h.s. of   \eqref{itocomp} can be rewritten as a jump term for the process $\nu^{i,k}$, namely
\begin{equation*}
\begin{split}
\bigg\langle \nu^{i,K}_{s}, 
 \big(r_{i} - \sum_{j=1}^M C^{ij}* \nu^{j,K}_{s}\big)f(s,\cdot)  \bigg\rangle = & K\int \nu^{i,K}_{s}(dx) \bigg[ 
  \bigg(\langle \nu^{i,K}_{s} +{1\over K}  \delta_x, f(s,\cdot) \rangle - \langle \nu^{i,K}_{s} , f(s,\cdot) \rangle  \bigg)r_{i}(x)\\
  & +  \bigg(\langle \nu^{i,K}_{s} - {1\over K} \delta_x , f(s,\cdot) \rangle - \langle \nu^{i,K}_{s} , f(s,\cdot) \rangle  \bigg)   \sum_{j=1}^M C^{ij}* \nu^{j,K}_{s}(x) \bigg] \\
 \end{split}
 \end{equation*}

It follows by It\^o's formula that for each function $F:{\cal  M}^M\to \RR$ of the form
\begin{equation}
  \label{cylind}
  F(\nu)=F(\langle \nu^1,f_{1}\rangle,\cdots, \langle \nu^M,f_{M}\rangle)
\end{equation}
with  $F\in {\cal C}^2_b(\mathbb{R}^M)$ and   $f_{i}\in
{\cal C}^{2}(\mathbb{R}^d)$  for $i\in \{1,\cdots,M\}$, the process
\begin{equation*} 
F( \nu_t^K) - F( \nu_0^K) -\int_0^t   L^KF (\nu_s^K) ds
\end{equation*} 
is a local martingale, if
\begin{align*}
  L^KF (\nu)&=L^K_e F(\nu)+L_d F(\nu)\label{generator}
\end{align*}
is the sum  of a jump operator $L^K_e F$ defined  by
\begin{align*}
  L^K_{e}F(\nu)&=K \sum_{i=1}^M\,\int \nu^{i}_{s}(dx) \,  r_{i}(x) \left(F(\nu+{1\over K}\delta_{x}^i)-F(\nu)\right)
 \notag\\
  &+K \sum_{i=1}^M \int \nu^{i}_{s}(dx) \,  (\sum_{j=1}^M C^{ij} *\nu^j(x))\left(F(\nu-{1\over K}\delta_{x}^i)-F(\nu)\right),
\end{align*}
where  the  $i$-th coordinate  of  the vector meaure $\delta_x^i\in {\cal M}^M$ is $\delta_x$ and the other ones vanish, plus a diffusion operator $L_d F$   given by
\begin{equation*}
  \label{gendiff}
  \begin{split}
   L_d F(\nu) =&\sum_{i=1}^M  \bigg(\big\langle\nu^i,  {\cal L}^i_{\nu_s}f_i 
 \big\rangle \partial_i \, F\, (\langle \nu^1,f_{1}\rangle,\cdots, \langle \nu^M,f_{M}\rangle) \\
  &+ \big\langle\nu^i,   ( \nabla f_{i} )^* a^i (\cdot,G^{i1}* \nu^i,\cdots,G^{iM}* \nu^M)  \nabla  f_{i}    \big\rangle
 \,  \partial^2_{ii}\,  F (\langle \nu^1,f_{1}\rangle,\cdots, \langle \nu^M,f_{M}\rangle)\big) \,\bigg)
 \end{split}
 \end{equation*}
(which is similar to the one obtained for branching diffusing processes  in \cite{D}.)

The  cylindrical functions \eqref{cylind} are a standard class generating the set of
bounded and measurable functions from $({\cal M}_{K})^M$ into $\mathbb{R}$. It is classical to deduce from the previous  that  the right continuous process   $(\nu^K_{t})_{t}$ in  $({\cal M}_{K})^M$ is  Markov, with law  characterized by the infinitesimal generator $L^K$.

\section{Large population approximation and non local Lotka-Volterra cross diffusion system   }
 
\subsection{Existence and uniqueness of weak measure solutions}

We now state a large population approximation for the previous $M$ species model by making the charge capacity $K$ tend to infinity. This result  in particular implies the  existence of weak solutions to a non local cross-diffusion system of nonlinear partial differential equations.

\begin{theo}\label{existunique}
Assume that for some $p\geq 3$,  $\sup_{K} \mathbb{E}(\langle \nu_{0}^{i,K},1\rangle^p)<+\infty $ for any $i=1, \cdots, M$. Assume also  $(\mathbf{H})$  and moreover that the sequence  of finite measures $( \nu^{1,K}_0 ,\cdots, \nu^{M,K}_0)$ converges in law  as $K$ goes to infinity  to the deterministic vector of finite measures $(\xi^1_0, \cdots, \xi^M_0)$. Then, the sequence  $( \nu^{1,K}_. ,\cdots, \nu^{M,K}_.)$ converges in law in  $\mathbb{D}([0,T],{\cal M}^M)$ when $K$ tends to infinity to the unique deterministic continuous finite measure-valued function $\xi=(\xi^1_{\cdot}, \cdots, \xi^M_{\cdot})$ which is a weak solution of the following cross-diffusion system: for  each $i=1, \cdots, M$, 
\begin{equation}\label{cross-diff}
\langle  \xi^{i}_{t},  f(t,\cdot)\rangle=   \langle \xi^i_{0}, f(0,\cdot) \rangle  + 
 \int_{0}^t 
\bigg\langle \xi^{i}_{s}, \partial_s f (s,\cdot)+ {\cal L}^i_{\xi_s}f (s,\cdot) + \big(r_{i}- \sum_{j=1}^M C^{ij}* \xi^{j}_{s}\big)f (s,\cdot) \bigg\rangle ds 
\end{equation}

for  every   $f\in{\cal C}_b^{1,2}([0,T]\times \mathbb{R}^d)$  such that
\begin{equation}\label{growthf}
\sup_{(t,x)\in [0,T]\times \RR^d}(1+ |x|) |\nabla f(t,x)|<\infty.
\end{equation}
\end{theo}

\medskip

\begin{rem}\label{examples} In contrast to the models of cross-diffusion  introduced by  \citeasnoun{SKT} or \citeasnoun{MK},
Equation \eqref{cross-diff} allows for long range interaction in the  coefficients of spatial diffusion. For example, taking $G^{ij}=H^{ij}=1$, the spatial behavior of individuals of  type $i$ depend on the total mass of the  subspecies $j$.   It also covers some cases where the diffusion matrix might vanish , e.g. $a^i(x,v_1,\dots,v_M)=I_d\Psi_i^2(\sum_{j=1}^M v_j)$ with $\Psi_i:[0,\infty]\to \RR_+$ a Lipschitz continuous function vanishing at $0$ and $I_d$ the identity matrix.

\end{rem}

\medskip

The proof of Theorem \ref{existunique}  consists in the following steps: 
\begin{itemize}
\item[i)] Propagation of moments of the total mass and control of  tails. 
\item[ii)] Uniform tightness of the  laws of $( \nu^{1,K}_. ,\cdots, \nu^{M,K}_.)$.
\item[iii)]  Identification of the limits in distribution as  solutions to \eqref{cross-diff}.
\item[iv)] Uniqueness of  solutions  to \eqref{cross-diff}. 
\end{itemize}

The  first  three steps are relatively standard and can be done by  general arguments developed in \citeasnoun{FM04},  \citeasnoun{CM07}, \cite{JourMelWo} or in \cite[Theorem 7.4]{BM}. In the remainder of this section, we provide  complete detailed proof  of steps i), ii) and iii) for the specific model we deal with. Step iv), namely  uniqueness   of weak measure solutions to  \eqref{cross-diff},  will  in turn require new  techniques and arguments, in order to deal with  the  nonlinear  diffusion terms of the equation.  The proof will be given in Section \ref{Uniqueness}.

\medskip

 The following bounds following from assumption  $(\mathbf{H})$ will be useful at several points: 
 for all $x\in \RR^d$ and   $ v_{j}\in \RR_+ , j\in \{1,\cdots, M\}$, we have
\begin{equation}\label{boundcoeftight}
|a^i(x,v_{1},\cdots, v_{M})| \leq C \Big(1+\sum_{j=1}^M v_{j}^2\Big)\, \mbox{ 
and  } \, |b^i(x,v_{1},\cdots, v_{M})| \leq C \Big(1+ |x|+\sum_{j=1}^M v_{j}\Big). 
\end{equation}
It follows that for any  $\nu=(\nu^i)_{1\leq i\leq M} \in {\cal M}^M$ and any function $f\in {\cal C}^2(\RR^d)$ satisfying the growth estimate \eqref{growthf},  one has
 \begin{equation}\label{boundtightL}
 \Big| {\cal L}^i_{\nu}f(x)\Big| \leq   C_f \Big(1+\sum_{j=1}^M\langle \nu^{j}, 1\rangle^2  \Big)
\end{equation}
and
 \begin{equation}\label{boundtightquadvarB}
   ( \nabla f)^* a^i (\cdot,G^{i1}* \nu^{i},\cdots,G^{iM}* \nu^{M})  \nabla f  \leq C_f \Big(1+\sum_{j=1}^M\langle \nu^{j}, 1\rangle^2  \Big)
\end{equation}
 for some finite constant $C_f>0$ depending only on the  supremum \eqref{growthf}.

   In the sequel, $C>0$ denotes some constant, which can change from line to line and which only  depends  on the assumptions of our model. 
  
\subsection{Moments of the total mass and  control of tails}

We start with two basic but crucial estimations. 

\begin{lem}\label{momentsmass}
Suppose $(\mathbf{H})$ holds.  Then, for each $T>0$, any $p\geq 1$  and  some constant $C_{T,p}>0$ we have
 $$\sup_{K} \mathbb{E}\left(\sup_{t\in [0,T]} \langle \nu_{t}^{i,K},1\rangle^p\right)\leq C_{T,p} \sup_{K} \mathbb{E}(\langle \nu_{0}^{i,K},1\rangle^p).  $$
for all $i=1, \cdots, M$.
\end{lem}
{\bf Proof.}  We can assume that $\sup_{K} \mathbb{E}(\langle \nu_{0}^{i,K},1\rangle^p)<\infty$.  Since the process $\langle \nu_{t}^{i,K},1\rangle^p$ is pure jump and of finite variation,  for all $t\geq 0$ we have
\begin{equation*}
\langle \nu_{t}^{i,K},1\rangle^p =\langle \nu_0^{i,K},1\rangle^p+ \sum_{s\leq t} \left[ \langle \nu^{i,K}_{s-},1\rangle +\Delta_s \right]^p- \langle \nu^{i,K}_{s-} , 1 \rangle^p ,
\end{equation*}
where $\Delta_s = {1\over K} $ or $- {1\over K} $ if $s$ is a jump time and  $\Delta_s=0$ otherwise.   Let $S_R$ be  the same stopping time as  in \eqref{mass}. Neglecting the negative jumps and  bounding the rate of positive ones,   we get that 
\begin{equation*}
\begin{split}
 \langle \nu^{i,K}_{t\wedge S_R}, 1\rangle^p \leq &\,   \langle \nu^{i,K}_{0}, 1\rangle^p + 
  \sum_{s\leq t \wedge S_R, \Delta_s={1\over K}} \left[ \langle \nu^{i,K}_{s-},1\rangle +{1\over K} \right]^p- \langle \nu^{i,K}_{s-} , 1 \rangle^p \\
  \leq  & \,  \langle \nu^{i,K}_{0}, 1\rangle^p \\
 &+ \int_{[0,t\wedge S_R]\times \mathbb{N}\times  \mathbb{R}_{+}}\left(\left[ \langle \nu^{i,K}_{s-},1\rangle +{1\over K} \right]^p- \langle \nu^{i,K}_{s-} , 1 \rangle^p \right)   {\bf 1}_{ \theta \leq \bar{r}, n\leq  K\,\langle \nu^{i,K}_{s-},1\rangle}  Q^i(ds,dn,d\theta).\\
 \end{split}
\end{equation*}
Since  for $a\geq 0$,  $b\in [0,1]$ and some constant $C_p>0$, 
$$(a+b)^p-a^p =\int_a^{a+b} p y^{p-1}dy \leq p \, b(a+b)^{p-1}\leq C_p \, b (a^{p-1}+1)$$
we deduce, in a similar way as in \eqref{bycompen},  that
\begin{equation*} 
\begin{split}
 \EE\left(\sup_{t\in [0,T\wedge S_R]} \langle \nu^{i,K}_{t}, 1\rangle^p\right)\leq   & \,  \EE\left( \langle \nu^{i,K}_{0},1\rangle^p \right)
 +\bar{r} C_p \int_0^T \EE\left( \sup_{s\in [0,t\wedge S_R]}  \langle \nu^{i,K}_{s} , 1 \rangle^p  + \langle \nu^{i,K}_{s} , 1 \rangle  \right)dt	\\
\leq  & \,  \EE\left( \langle \nu^{i,K}_{0},1\rangle^p \right)
 +\bar{r} C_p \int_0^T \EE\left( 1+2\sup_{s\in [0,t\wedge S_R]}  \langle \nu^{i,K}_{s} , 1 \rangle^p  \right)dt.	\\
\end{split}
\end{equation*}
The result follows from  Gronwall's  lemma, letting then $R\nearrow \infty$. 

\Box

\bigskip 

Let  $B(0,R)$ denote the centered ball of radius $R>0$ in $\RR^d $ and  $\Psi:\RR^d\to[0,1]$ be  a  radially  non-decreasing function of class $C^2(\RR^d)$ such that $\Psi(x)=0$ for all $x\in B(0,1)$ and  $\Psi(x)=1$ for all $x\in B(0,2)^c$.  For each $m\geq 1$  define a $C^2(\RR^d)$ function $\Psi_m$ by   
\begin{equation}\label{Psim}
\Psi_m(x):= \Psi(x/m), \quad x\in \RR^d.
\end{equation}
Observe that, since for each $m\geq 0$ and any  $\mu\in {\cal M}$ one has
$$\mu(B(0,2m)^c)  \leq  \langle \mu, \Psi_m\rangle \leq \mu(B(0,m)^c),$$
a family  $(\mu_k)_{k} $ of elements of ${\cal M}$ is tight if and only if 
 $\lim_{m\to \infty }\sup_{k} \langle \mu_k, \Psi_m\rangle =0$.

\begin{lem}\label{tightsup}
 Under the assumptions  of  Theorem \ref{existunique}, there is a constant $C_0 >0$ depending  on  $T$ and on $\sum_{i=1}^{M} 
 \sup_{K} \mathbb{E}(\langle \nu_{0}^{i,K},1\rangle^3)<\infty$ but not on 
 $m$,  such that
$$ \sup_{K} \mathbb{E}\left(\sup_{t\in [0,T]} \langle \nu_{t}^{i,K},\Psi_m\rangle\right)\leq  C_0  \sup_{K} \mathbb{E}\left( \langle \nu_0^{i,K},\Psi_m\rangle\right) $$
for all  $i=1, \cdots, M$. 
As a consequence, for  all  $i=1, \cdots, M$ we have
$$\lim_{m\to \infty} \sup_{K} \mathbb{E}\left(\sup_{t\in [0,T]} \langle \nu_{t}^{i,K},\Psi_m\rangle\right)=0.$$
\end{lem}

{\bf Proof:}   It is readily seen that the functions $(1+ |x|) \nabla \Psi_m(x)$  and $(1+ |x|^2) \operatorname{Hess} \Psi_m(x)$ are globally bounded, uniformly in $m\geq 1$.  
 Therefore,  \eqref{boundtightL} holds for some constant $C_{\Psi_m}$ which in fact does not depend on $m$. 
  This and \eqref{itocomp} yield
\begin{equation*}
\begin{split}
\langle  \nu^{i,K}_{t},  \Psi_m\rangle\leq    &\langle \nu^{i,K}_{0}, \Psi_m\rangle  + C
 \int_{0}^{t}  \langle   \nu^{i,K}_{s},    \Psi_m \rangle + \langle   \nu^{i,K}_{s},   1\rangle 
\left( 1 + \sum_{j=1}^M \langle \nu^{j,K}_{s},1 \rangle^2\right)  ds \\
& + M^{i,\Psi_m,B}_{t}+ M^{i,\Psi_m,Q}_{t} \,  ,\\
\end{split}
\end{equation*}
 where  the  local martingales $M^{i,\Psi_m,Q}_{t}$ and $ M^{i,\Psi_m,B}_{t}$  moreover satisfy 
  $$  \langle M^{i,\Psi_m,Q}\rangle_{t} +   \langle M^{i,\Psi_m,B}\rangle_{t}  \leq {C\over K}   \int_{0}^{t}  \langle   \nu^{i,K}_{s},   1\rangle 
\left( 1 + \sum_{j=1}^M \langle \nu^{j,K}_{s},1 \rangle  + \langle \nu^{j,K}_{s},1 \rangle^2 \right) ds$$
for some constant which is independent of $m$.
Using the BDG inequality,  we get  for all $i=1,\dots,M$ and $K>0$  that
\begin{equation*}
\begin{split}
 \,  \EE\left( \sup_{t\in [0,T]}  \langle \nu^{i,K}_{t},  \Psi_m\rangle\right)\leq    &  \,  \EE\left( \langle \nu^{i,K}_{0}, \Psi_m\rangle\right)  + C
 \int_{0}^{T}    \EE\left( \sup_{s\in [0,t]}  \langle   \nu^{i,K}_{t},    \Psi_m \rangle\right) dt \\
 & +  C T 
 \EE\left(\sup_{t\in [0,T]}   \langle   \nu^{i,K}_{t},   1\rangle 
\left( 1 +  \sum_{j=1}^M  \sup_{t\in [0,T]}  \langle \nu^{j,K}_{t},1 \rangle^2  \right)  \right)\\
 & +  C \sqrt{ T }
 \EE\left( \sqrt{\sup_{t\in [0,T]}   \langle   \nu^{i,K}_{t},   1\rangle 
\left( 1 +  \sum_{j=1}^M  \sup_{t\in [0,T]}  \langle \nu^{j,K}_{t},1 \rangle^2  \right)  }\right) . \\
\end{split}
\end{equation*}
We thus get
\begin{equation*}
\begin{split}
 \,  \EE\left( \sup_{t\in [0,T]}  \langle \nu^{i,K}_{t},  \Psi_m\rangle\right)
\leq    &  \,  \EE\left( \langle \nu^{i,K}_{0}, \Psi_m\rangle\right)  + C
 \int_{0}^{T}    \EE\left( \sup_{s\in [0,t]}  \langle   \nu^{i,K}_{t},    \Psi_m \rangle\right) dt \\
  & +  C (T + 1 )   \EE\left( 1+\sup_{t\in [0,T]}   \langle   \nu^{i,K}_{t},   1\rangle^3\right)^{1/3} 
   \EE\left( \, 1 + \sum_{j=1}^M  \sup_{t\in [0,T]} \langle \nu^{j,K}_{t},1 \rangle^3  \right)^{2/3} \\
   \leq    &   \EE\left( \langle \nu^{i,K}_{0}, \Psi_m\rangle\right)  + C (T + 1 )  
   \EE\left( 1 +  \sum_{j=1}^M    \langle \nu^{j,K}_{0},1 \rangle^3  \right) \\
 &   + C
 \int_{0}^{T}    \EE\left( \sup_{s\in [0,t]}  \langle   \nu^{i,K}_{t},    \Psi_m \rangle\right) dt, \\ 
\end{split}
\end{equation*}
where we used Lemma \ref{momentsmass} in the last inequality. The first statement follows with Gronwall's lemma.

As for the last statement, notice that  for each $n\in \NN^*$, by Lemma \ref{momentsmass}  with $p=3$ and the inequality  $z\leq z\wedge n+ z^3n^{-2} $ for $z\geq 0$ we get
$$\EE( \langle  \nu_{0}^{i,K},\Psi_m\rangle)\leq  \EE( \langle  \nu_{0}^{i,K},\Psi_m\rangle  \wedge n) + {C\over n^2}, $$
from where $$ \limsup_{K\to \infty} \mathbb{E}\left( \langle \nu_0^{i,K},\Psi_m\rangle\right) \leq   \langle  \xi_{0}^{i},\Psi_m\rangle.  $$
This easily yields  $\lim_{m\to \infty}\sup_{K\in \NN^*} \EE( \langle  \nu_{0}^{i,K},  \Psi_m \rangle)$ and we conclude using the first part.

\Box

\subsection{Tightness }
Let us denote by ${\cal P}_K$ the law of the process $( \nu^{1,K}_. ,\cdots, \nu^{M,K}_.)$.
We will  first prove a slightly weaker tightness result   than required. 

\begin{lem}\label{tightvague}
Let $\tilde{{\cal P}}_K$  denote the law of the process $( \nu^{1,K}_. ,\cdots, \nu^{M,K}_.)$, $K\in\NN^*$, when  seen  as a random element of the space   $(\mathbb{D}([0,T],{\cal M}_v))^M$, with ${\cal M}_v$ denoting the space $ {\cal M}$ endowed with the vague topology.  Then, under the assumptions  of  Theorem \ref{existunique},  the sequence  $(\tilde{{\cal P}}_K)_{K\in \NN^*}$ is tight. 
\end{lem}

Since the product of compact sets is a compact set, it is enough to prove  tightness  in $\mathbb{D}([0,T],{\cal M}_v)$  for the  family of processes  $\nu^{i,K}_. $  , $K\in\NN^*$, for each $i=1,\dots,M$. 
To that end, we  use a  criterion established in \cite{Roell}:

\medskip

Let  $(\mu^K_t)_{t\in [0,T]}$  denote the canonical process in  $(\mathbb{D}([0,T],{\cal M}_v))$ under a given probability law ${\cal Q}_K$.  A sequence  $({\cal Q}_K)_{K\in \NN^*}$ of  laws  on  $(\mathbb{D}([0,T],{\cal M}_v))$ is tight  if and only if
\begin{itemize}
\item    For each function $f\in {\cal C}_0^2(\RR^d)$ the sequence of laws   $({\cal Q}^f_K)_{K\in \NN^*}$ of the processes $(\langle  \mu^K_{t},  f\rangle)_{t\in [0,T]}$  is tight in  $\mathbb{D}([0,T],\RR)$.
\item The sequence of laws   $({\cal Q}^1_K)_{K\in \NN^*}$ of the processes $(\langle  \mu^K_{t},  1\rangle)_{t\in [0,T]}$  is tight  in  $\mathbb{D}([0,T],\RR)$.
\end{itemize}
\medskip

{\bf Proof of Lemma \ref{tightvague}:}  For  each $i=1,\dots, M$, we need to prove that the  sequence of  laws of the processes $(\langle  \nu^{i,K}_{t},  f\rangle)_{t\in [0,T]}$,$K\in \NN^*$,  is tight  in  $\mathbb{D}([0,T],\RR)$  for all $f\in {\cal C}_0^2(\RR^d)$  and  for $f=1$.   The same proof will work in both cases.  Let us write
$$A^f_t:= \int_{0}^t 
\bigg\langle \nu^{i,K}_{s},{\cal L}^i_{\nu_s^K}f  + \big(r_{i}- \sum_{j=1}^M C^{ij}* \nu^{j,K}_{s}\big)f  \bigg\rangle ds.$$ 
By \eqref{itocomp} and the Aldous-Rebolledo criterion, it is enough to prove that  
$$ \lim_{\eta\to 0}\limsup_{K\to \infty} \PP\left(\sup_{t\in [0,T]} \langle \nu_{t}^{i,K},f\rangle>\frac{1}{\eta}\right)=0, $$
and that, if $S_\delta$ and $S'_\delta$  denote  for each $\delta>0$  generic $({\cal F}_t)$-stopping times such  that $S_\delta\leq S'_\delta\leq (S_\delta+\delta)\wedge T$,  then
 $$  \lim_{\delta\to 0} \limsup_{K\to \infty} \sup_{S_\delta,S'_\delta}\PP(  |A^f_{S_\delta}-A^f_{S'_\delta}|\geq \eta)=0$$
and
$$  \lim_{\delta\to 0} \limsup_{K\to \infty} \sup_{S_\delta,S'_\delta}\PP(  | \langle M^{i,f,B}\rangle_{S_\delta}- \langle M^{i,f,B}\rangle_{S'}|+  | \langle M^{i,f,Q}\rangle_{S_\delta}- \langle M^{i,f,Q}\rangle_{S'_\delta}|\geq \eta)=0.$$

The first property follows form Lemma \ref{momentsmass} and Markov's inequality. From \eqref{boundtightL} and using similar estimates as in the proof of Lemma \ref{tightsup},  we then get 
\begin{equation*}
\begin{split}
\EE\left(  |A^f_{S_\delta}-A^f_{S'_\delta}|\right) \leq    &  C  \EE\left(
 \int_{S_\delta}^{S'_\delta}  \langle   \nu^{i,K}_{s},   1\rangle 
\left( 1 + \sum_{j=1}^M \langle \nu^{j,K}_{s},1 \rangle^2\right)  ds\right) \\
   \leq     &  \delta C
   \EE\left( 1 +  \sum_{j=1}^M  \sup_{t\in [0,T]}  \langle \nu^{j,K}_{t},1 \rangle^3  \right). \\
 \end{split}
\end{equation*}
Moreover,  \eqref{quadvarB}, \eqref{boundtightquadvarB} and  \eqref{quadvarQ}  similarly yield
\begin{equation*}
\begin{split}
\EE\bigg(   | \langle M^{i,f,B}\rangle_{S_\delta}- \langle M^{f,B}\rangle_{S'_\delta}|+  | \langle M^{f,Q} & \rangle_{S_\delta}-  \langle M^{i,f,Q}\rangle_{S'_\delta}|\bigg)    \\  
\leq &  C \, \sqrt{\delta} \EE\left( \sqrt{\sup_{t\in [0,T]}   \langle   \nu^{i,K}_{t},   1\rangle 
\left( 1 +  \sum_{j=1}^M  \sup_{t\in [0,T]}  \langle \nu^{j,K}_{t},1 \rangle^2  \right)  }\right) \\
 \leq     & \sqrt{ \delta} C
   \EE\left( 1 +  \sum_{j=1}^M  \sup_{t\in [0,T]}  \langle \nu^{j,K}_{t},1 \rangle^3  \right) \\
 \end{split}
\end{equation*}
The remaining conditions  then follow  by Markov's inequality. 

\Box

\bigskip

Tightness in  $\mathbb{D}([0,T],{\cal M}^M)$ can now be deduced  from   the next criterion  proved in \cite{MelRoell}:  

\medskip
  
For each $j\in \NN$,   let  $(\mu^j_t)_{t\in [0,T]}$  denote the canonical process in  $\mathbb{D}([0,T],{\cal M})$ under a given probability law ${\cal Q}_j$ (notice that  $t\mapsto \mu^j_t$ is assumed  to be weakly continuous under ${\cal Q}_j$), and ssume that
\begin{itemize}
\item    the sequence  $({\cal Q}_j)_{j\in \NN}$  is weakly convergent in $\mathbb{D}([0,T],{\cal M}_v)$ to some  law ${\cal Q}$, 
\item the canonical process $(\mu_t)_{t\in [0,T]}$  in  $\mathbb{D}([0,T],{\cal M})$ is continuous under ${\cal Q}$ and
\item   the sequence of laws   of the processes $(\langle  \mu^j_{t},  1\rangle)_{t\in [0,T]}$ weakly converge to  the law   of the process $(\langle  \mu_{t},  1\rangle)_{t\in [0,T]}$.
\end{itemize}
Then,  the sequence  $({\cal Q}_j)_{j\in \NN}$  is weakly convergent in $\mathbb{D}([0,T],{\cal M})$ to  ${\cal Q}$. 

 This criterion's assumptions hold in our case, thanks to the next result: 

\begin{lem}\label{tightweak}   Let  the law  $\tilde{{\cal P}}$  on $(\mathbb{D}([0,T],{\cal M}_v))^M$ be a  weak limit  of some subsequence  $(\tilde{{\cal P}}_{K_j})_{j\in \NN^*}$ and let $( \nu^{1}_. ,\cdots, \nu^{M}_.)$ denote the canonical process under   $\tilde{{\cal P}}$.  For each $i=1,\dots, M$ we have: 
\begin{itemize}
\item[i)] $ \mathbb{E}\left(\sup_{t\in [0,T]} \langle \nu_{t}^{i},1\rangle^3\right)<\infty $  and $\lim_{m\to \infty}\mathbb{E}\left(\sup_{t\in [0,T]} \langle \nu_{t}^{i},\Psi_m\rangle\right)=0$,
\item[ii)]  the process $t\mapsto \nu^i_t \in {\cal M}$ is a.s. weakly continuous under $\tilde{{\cal P}}$ and
\item[iii)] the processes $\langle \nu^{i,K_j}_., 1\rangle$ converge in distribution to  $\langle \nu^{i}_., 1\rangle$ in $\mathbb{D}([0,T],\RR)$.
\end{itemize}

\end{lem}

{\bf Proof:}  For notational simplicity we assume that  the whole sequence $( \nu^{1,K}_. ,\cdots, \nu^{M,K}_.)$  in  $(\mathbb{D}([0,T],{\cal M}_v))^M$ converges in distribution  to $( \nu^{1}_. ,\cdots, \nu^{M}_.)$. 

Recall that for each $\psi\in {\cal C}_0(\RR^d)$  the mappings 
$$\mu_{\cdot}\in \mathbb{D}([0,T],{\cal M}_v) \mapsto  \sup_{t\in [0,T]} |  \langle \mu_{t}, \psi \rangle -  \langle \mu_{t-}, \psi \rangle| \in \RR$$ 
and
$$\mu_{\cdot}\in \mathbb{D}([0,T],{\cal M}_v) \mapsto  \sup_{t\in [0,T]}   \langle \mu_{t}, \psi \rangle
\in \RR$$ 
are continuous, by standard properties of the Skorokhod topology. 
\medskip

i)  Fix $\phi\in {\cal C}_b(\RR^d)$ with values in $[0,1]$ and for  $l\in \NN$, set $\phi_l:=\phi(1-\Psi_l)\in  {\cal C}_0(\RR^d)$. From  Lemma \ref{momentsmass}, we get  that
$$   \mathbb{E}\left(n \wedge \sup_{t\in [0,T]} \langle \nu_{t}^{i},\phi_l\rangle^3\right)< C$$
for some constant  $C<\infty$ neither depending on $n$ nor on $l$.  The first property then  follows by   choosing  $\phi=1$ and
letting $n,l\to \infty$  using Fatou's lemma.  Taking next $\phi=\Psi_m$ for a fixed $m\geq 1$, we get by   Lemma  \ref{tightsup}   that  
$$ \mathbb{E}\left(n\wedge  \sup_{t\in [0,T]} \langle \nu_{t}^{i},\phi_l\rangle\right)  \leq 
 \sup_{K} \mathbb{E}\left(\sup_{t\in [0,T]} \langle \nu_{t}^{i,K},\Psi_m\rangle\right) . $$
Letting $n,l\to \infty$, thanks to Lemma \ref{tightsup} the second  asserted property follows.

\medskip

ii)  For each $\psi\in {\cal C}_0(\RR^d)$ we have 
\begin{equation*}
\sup_{t\in [0,T]} |  \langle \nu^{i}_t, \psi \rangle- \langle \nu^{i}_{t-}, \psi \rangle|\leq 
n \wedge  \sup_{t\in [0,T]} |  \langle \nu^{i}_t, \psi \rangle- \langle \nu^{i}_{t-}, \psi \rangle|  +\frac{2\|\psi\|_{\infty}}{n^2}\sup_{t\in [0,T]} |  \langle \nu^{i}_t, 1\rangle^3| \\
\end{equation*}
Since for all  $n,K\in \NN^*$    one has $\sup_{t\in [0,T]} |  \langle \nu^{i,K}_t, \psi \rangle- \langle \nu^{i,K}_{t-}, \psi \rangle|\wedge n \leq   \frac{\|\psi\|_{\infty}}{K}$, the  first term on the right hand side  is null. By part i),  the  last  term is a.s. finite and goes to $0$ when $n$ goes to infinity. Thus, $\sup_{t\in [0,T]} |  \langle \nu^{i}_t, \psi \rangle- \langle \nu^{i}_{t-}, \psi \rangle|=0$ a.s. and the statement follows by separability of $ {\cal C}_0(\RR^d)$. 

\medskip

iii)    Let  $F: \mathbb{D}([0,T],\RR)\to \RR$ be a Lipschitz bounded function and $i=1,\dots,M$  be fixed. It is enough to show that
$$\EE( F(\langle \nu^{i,K}_., 1\rangle))\to \EE(F(\langle \nu^{i}_., 1\rangle))$$
when $K\to \infty$. We have for each $K,m\geq 1$, 
\begin{equation*}
\begin{split}
|\EE( F(\langle \nu^{i,K}_., 1\rangle))- \EE (F(\langle \nu^{i}_., 1\rangle))|\leq & \EE | F(\langle \nu^{i,K}_., 1\rangle) - F(\langle \nu^{i,K}_., 1-\Psi_m\rangle)| \\
& +  |\EE( F(\langle \nu^{i,K}_., 1-\Psi_m \rangle))-  \EE(F(\langle \nu^{i}_., 1-\Psi_m\rangle))|\\
& +  \EE | F(\langle \nu^{i}_., 1\rangle-  F(\langle \nu^{i}_., 1-\Psi_m\rangle)|\\
\leq & C\sup_{K'} \EE\left(\sup_{t\in [0,T]} \langle \nu_{t}^{i,K'},\Psi_m\rangle\right)\\
& +  |\EE( F(\langle \nu^{i,K}_., 1-\Psi_m \rangle))-  \EE(F(\langle \nu^{i}_., 1-\Psi_m\rangle))|\\
& + C \EE\left(\sup_{t\in [0,T]} \langle \nu_{t}^{i},\Psi_m\rangle\right),\\
\end{split}
\end{equation*}
where, in the second inequality, we used the fact that the Skorokhod distance in  $\mathbb{D}([0,T],{\cal M})$ is bounded above by the distance induced by the supremum norm. By Lemma \ref{tightsup}  and part i),  for each $\varepsilon>0$  there is $m\geq 0$ large enough such that the first and third terms on the right hand side of the last inequality are less that $\varepsilon$.  On the other hand, the  expectation of the second term goes to $0$ when $K\to \infty$ since $1-\Psi_m\in {\cal C}_0(\RR^d)$ and the mapping $\mu \in \mathbb{D}([0,T],{\cal M}_v)\mapsto \langle \mu_., 1-\Psi_m\rangle \in \mathbb{D}([0,T],\RR)$ is continuous.  We thus get
\begin{equation*}
\limsup_{K\to \infty} |\EE( F(\langle \nu^{i,K}_., 1\rangle))- \EE (F(\langle \nu^{i}_., 1\rangle))|\leq 2\varepsilon
\end{equation*}
and the statement follows.

\Box

Bringing all together, we get:
\begin{prop}
Under the assumptions  of  Theorem \ref{existunique},  the sequence  $({\cal P}_K)_{K\in \NN^*}$ of laws on  $\mathbb{D}([0,T],{\cal M}^M)$ is tight.
\end{prop} 

 {\bf Proof:}   By Lemma  \ref{tightvague}
 and Lemma \ref{tightweak} applied to a subsequence, we readily get that   the sequence  $({\cal P}_K)_{K\in \NN^*}$ is tight when each process  $( \nu^{1,K}_. ,\cdots, \nu^{M,K}_.)$ is seen  as a random element of the product space   $(\mathbb{D}([0,T],{\cal M}))^M$. Taking as distance in ${\cal M}^M$ the maximum of the distances used in each copy of ${\cal M}$, the Skorokhod modulus of continuity of an element $\mu\in \mathbb{D}([0,T],{\cal M}^M)$ is  seen to be bounded above by the maximum of  the $M$   modulus of continuity of its coordinates $\mu^i\in \mathbb{D}([0,T],{\cal M})$.
 The result easily follows.
 \Box
 
\subsection{Identification of  the limiting processes}\label{identlimitweak}

\begin{lem}  Under the assumptions  of  Theorem \ref{existunique}, let  the law  ${\cal P}$ be a  weak limit   on $\mathbb{D}([0,T],{\cal M})^M$   of some subsequence  $({\cal P}_{K_j})_{j\in \NN^*}$ and let $( \nu^{1}_. ,\cdots, \nu^{M}_.)$ denote the canonical process under  ${\cal P}$.  Then,  ${\cal P-}$a.s.  we have 
\begin{equation}\label{eqident}
\langle  \nu^{i}_{t},  f\rangle=   \langle \xi^i_{0}, f\rangle  + 
 \int_{0}^t 
\bigg\langle \nu^{i}_{s}, \partial_s f (s,\cdot)+ {\cal L}^i_{\nu_s}f (s,\cdot) + \big(r_{i}- \sum_{j=1}^M C^{ij}* \nu^{j}_{s}\big)f (s,\cdot) \bigg\rangle ds 
\end{equation}
for all $i=1,\dots,M$, every $t\in [0,T]$ and all $f\in {\cal C}_b^{1,2}([0,T], \RR^d)$  such that 
$\sup_{(t,x)\in [0,T]\times \RR^d}(1+ |x|^2)\left( |\nabla f(t,x)|^2+|\operatorname{Hess} f(t,x)|\right)<\infty$ and every $t\in [0,T]$. 
\end{lem}

{\bf Proof:}   For notational simplicity we again assume that the whole sequence   $({\cal P}_{K})_{K\in \NN^*}$   converges to ${\cal P}$. 
For fixed  $i=1,\dots,M$, $t\in [0,T]$   and $f\in {\cal C}_b^{1,2}([0,T], \RR^d)$ with the required growth conditions, define a functional   $\Phi^i: \mathbb{D}([0,T],{\cal M}^M)\to \RR_+$ acting on an element $\mu=( \mu^{1}_. ,\cdots, \mu^{M}_.)$  by
\begin{equation*}
\begin{split}
\Phi^i(\mu^{1}_. ,\cdots, \mu^{M}_.)=  & \langle  \mu^{i}_{t},  f\rangle -  \langle \xi^i_{0}, f\rangle \\
& - 
 \int_{0}^t 
\bigg\langle \mu^{i}_{s}, \partial_s f (s,\cdot)+ {\cal L}^i_{\mu_s}f (s,\cdot) + \big(r_{i}- \sum_{j=1}^M C^{ij}* \mu^{j}_{s}\big)f (s,\cdot) \bigg\rangle ds \\
\end{split}
\end{equation*}
 Since, from the previous section,   ${\cal P}$ is concentrated on $ \mathbb{C}([0,T],{\cal M}^M) $, it is enough to show that for each $n\in \NN$, the  bounded function  $\mu\mapsto |\Phi^i(\mu) |\wedge n$  is continuous at points $\mu\in  \mathbb{C}([0,T],{\cal M}^M) $ and  that
 $\EE(|\Phi^i(\nu^{1,K}_. ,\cdots, \nu^{M,K}_.))|\wedge n)\to 0$ when $K\to \infty$. 
 
 Let $(\mu_{\cdot}(k))_{k\in \NN}$ be a sequence converging in $\mathbb{D}([0,T],{\cal M}^M)$ to $\mu_{\cdot}\in  \mathbb{C}([0,T],{\cal M}^M) $.   Then, $\mu^i_{\cdot}(k)\to \mu^i_{\cdot}$ for each $i=1,\dots, M$, uniformly on $[0,T]$ when $k\to \infty$ and,   in particular, all masses are bounded uniformly in time and  in $k$.  Since for each $t\in [0,T]$ the coordinates of the sequence $(\mu_{t}(k))_{k\in \NN}$ are tight, the sequences  $(\mu_{t}^i(k)\otimes \mu_{t}^j(k) )_{k\in \NN}$ and $(\mu_{t}^i(k)\otimes \mu_{t}^j(k) \otimes \mu_{t}^j(k)  )_{k\in \NN}$  are tight too, and by  taking test functions of product form they are seen to converge weakly  to the product measures $\mu_{t}^i\otimes \mu_{t}^j $ and $\mu_{t}^i\otimes \mu_{t}^j \otimes \mu_{t}^j  $.    These facts and the (locally) Lipschitz character of the functions $a^i$ and $b^i$ easily yield the continuity of $ |\Phi^i |\wedge n$ at $\mu$, by dominated convergence.  

Last, from  equations \eqref{itocomp},   \eqref{quadvarB} and  \eqref{quadvarQ},  the bound \eqref{boundtightquadvarB} and similar estimates as in the last part of the proof of Lemma \ref{tightvague}, we get that 
$$\EE(|\Phi^i(\nu^{1,K}_. ,\cdots, \nu^{M,K}_.))|\wedge n)\leq {C\over\sqrt{ K}}$$ for all $n\in \NN$. This together with the convergence in law of the initial conditions complete the proof. 

\Box

\begin{rem}
\label{bound}
Notice for later reference that any  solution to \eqref{cross-diff} satisfies
$\sup_{t\in [0,T]} \|\xi_t^i\|_{TV}\leq e^{{\bar r}_i T}  \|\xi_0^i\|_{TV} $ for $i=1,\dots,M$, as is readily seen by taking $f=1$ and using the non negativity of the functions $C^{ij}$ and Gronwall's lemma.
\end{rem}

\section{Uniqueness}\label{Uniqueness}

  {Uniqueness is established in next result where, for notational simplicity, we will deal  only with the case $M=2$. All arguments easily extend to the general case.

\medskip

\begin{prop}\label{uniqueness}
Let  $\xi=(\xi^1,\xi^2)$ and
$\tilde{\xi}=(\tilde{\xi^1},\tilde{\xi^2})$ be two solutions of
the system \eqref{cross-diff} in $[0,T]$ with $M=2$. Then
$(\xi^1_t,\xi^2_t)=(\tilde{\xi}^1_t,\tilde{\xi}^2_t)$ for all $t
\in [0,T]$.
\end{prop}

We will  need to  use a distance that is  weaker than to total variation one
but better adapted  to perturbations in the diffusion coefficients. Denote by ${\cal L}{\cal B}(\mathbb{R}^d)$ the space of Lipschitz
continuous and bounded functions on $\mathbb{R}^d$, and by
$\|\cdot\|_{{\cal L}{\cal B}}$ or simply $\|\cdot\|$ the
corresponding norm,
$$\|\varphi\|_{{\cal L}{\cal B}} := \sup_{x\not = y}
\frac{|\varphi(x)-\varphi(y)|}{|x-y|}\quad +\quad  \sup_{x}
|\varphi(x)|.$$

We endow $ {\cal M}$ with the distance induced by the dual norm
with respect to ${\cal L}{\cal B}(\mathbb{R}^d)$,
$$\|\eta-\mu\|_{{\cal L}{\cal B}^*} := \sup_{\|\varphi\|_{{\cal L}{\cal
B}}\leq 1}| \langle\mu- \eta, \varphi\rangle|, \quad \eta,\mu \in {\cal
M}.$$

Given a solution $(\xi^1_t,\dots, \xi^M_t)_{t\in [0,T]}$ to \eqref{cross-diff}, let us set for $i=1,\dots, M$,
\begin{eqnarray}\label{coefflow}
\sigma(i,t,x)&=&\sigma^i(x,G^{i1}* \xi_t^1(x),\cdots, G^{iM}* \xi_t^M(x)) \, ,\notag\\
b(i,t,x)&=&b^i(x,H^{i1}* \xi_t^1(x),\cdots, H^{iM}* \xi_t^M(x)). 
\end{eqnarray}

\begin{rem}\label{coefflowlip} From assumption $(\mathbf{H})$ i),ii)  and Remark \ref{bound}, the functions in \eqref{coefflow} are Lipschitz functions of $x \in \RR^d$, uniformly in $[0,T]$.
\end{rem}

We introduce next the family of SDEs  associated with  the coefficients \eqref{coefflow}, and the corresponding transition semigroups. For each $x\in\RR^d$ and $s\in [0,T]$  consider   the unique (strong) solution $$X^i_{s,t}(x)=(X^{i,1}_{s,t}(x),\dots, X^{i,d}_{s,t}(x)) \quad t\in [s,T]$$ of
the stochastic differential equation in $\RR^d$, 
\begin{equation}\label{stochflow}
X^i_{s,t}(x)= x+\int_s^t \sigma(i,r,X^i_{s,r}(x)) dB_r^i + \int_s^t b(i,r,X^i_{s,r}(x)) dr ,\quad t\in [s,T]
\end{equation}
where $B^i=(B^{i,q})_{q=1}^d$ is a standard $d$-dimensional Brownian motion in a given probability space. 
   The fact that  for each $s$ the mapping  $(t,x)\mapsto X^i_{s,t}(x)$  is measurable can be classically deduced from the properties of  functions $\sigma(i,t,x)$ and $b(i,t,x)$ noted in Remark \ref{coefflowlip}.
 The three parameter process $(s,t,x)\mapsto X^i_{s,t}(x)$  is called the {\it stochastic flow} associated with the coefficients $\sigma(i,t,x)$ and $b(i,t,x)$. Finer properties  of this process will be recalled  and used later.

Given a second  solution   $(\tilde{\xi}^1_t,\tilde{\xi}^2_t)_{t\in [0,T]}$ of \eqref{cross-diff}, define analogously  coefficients $\tilde{\sigma}(i,t,x)$ and $\tilde{b}(i,t,x)$ in terms of  $(\tilde{\xi^1},\tilde{\xi^2})$, and the processes $\tilde{X}^i_{s,t}(x)$  given for $i=1,2$ by the  solution to the  SDEs  
\begin{equation*}
\tilde{X}^i_{s,t}(x)= x+\int_s^t \tilde{\sigma}(i,r,\tilde{X^i}_{s,r}(x)) dB_r^i+ \int_s^t \tilde{b}(i,r,\tilde{X^i}_{s,r}(x)) dr,
\end{equation*}
 driven by {\it the same} Brownian motions $B^i$ as  the processes $X^i_{s,t}(x)$ in \eqref{stochflow}. 

\medskip

The proof  of Proposition \ref{uniqueness}   will   rely on   stability  properties of  the non homogenous transition
 semigroups of $X^i_{s,t}(x) $  and $\tilde{X}^i_{s,t}(x)$, which we  respectively denote   by $$P_{s,t}^i(x,dy) =\PP( X^i_{s,t}(x)    \in dy) \mbox{  and }\tilde{P}^i_{s,t}(x,dy)=\PP( \tilde{X}^i_{s,t}(x)   \in dy).$$
 As usual,  $C>0$ denotes a constant that may change from line to line.

\begin{lem}\label{bornlip}
For all $T>0$, $\varphi \in {\cal L}{\cal
 B}(\RR^d)$, $0\leq s\leq t\leq T$, $i=1,2$, 
\begin{itemize}
\item[a)] $\|P_{s,t}^i\varphi \|^2_{{\cal L}{\cal B}}\leq C(\xi^1,\xi^2,T) \|\varphi \|^2_{{\cal L}{\cal B}}$
and $\|\tilde{P}_{s,t}^i\varphi \|^2_{{\cal L}{\cal B}}   \leq
C(\tilde{\xi^1},\tilde{\xi^2},T) \|\varphi \|^2_{{\cal L}{\cal B}}$ .
\item[b)]  For all $x \in \RR^d$,
\begin{equation*}
|P_{s,t}^i\varphi(x)-\tilde{P}^i_{s,t}\varphi(x)|^2\leq
C'(\xi^1,\xi^2,\tilde{\xi}^1,\tilde{\xi}^2 ,T)\|\varphi \|^2_{{\cal L}{\cal
B}}\displaystyle{\int_s^t \|\xi^1_r-\tilde{\xi}^1_r\|^2_{{\cal
L}{\cal B}^*}+ \|\xi^2_r-\tilde{\xi}^2_r\|^2_{{\cal
L}{\cal B}^*}dr} 
\end{equation*}
\end{itemize}
The constants depend on $\xi^i$  or  $\tilde{\xi}^i$ only
through $\sup_{t\in [0,T]} \|\xi_t^i\|_{TV}$ and $\sup_{t\in
[0,T]} \|\tilde{\xi}_t^i\|_{TV}$, and can be chosen to depend   only on $e^{{\bar r}_i T}  \|\xi_0^i\|_{TV} $. 
\end{lem}

{\bf Proof:} {\it a)} Is is enough to control the  Lipschitz
constants of the function $P_{s,t}^i\varphi$ or
$\tilde{P}_{s,t}^i\varphi$. We have, by Burkholder-Davis-Gundy inequality
\begin{equation*}
\begin{split}
\EE|X^i_{s,t}(x)-X^i_{s,t}(y)|^2\leq  & C|x-y|^2+C\int_s^t \EE|\sigma(i,r,X^i_{s,r}(x))-\sigma(i,r,X^i_{s,r}(y))|^2 dr \\
& \ +C\int_s^t \EE|b(i,r,X^i_{s,r}(x))-b(i,r,X^i_{s,r}(y))|^2 dr \\
 \leq & C|x-y|^2 +C\int_s^t
\EE|X^i_{s,r}(x)-X^i_{s,r}(y)|^2
dr\\
\end{split}
\end{equation*}
for all $s\leq t$. The above constants depend on bounds for  the Lipschitz constants of the coefficients $\sigma^i$, $b^i$ and   of the Kernels $G^{1j}$ and $H^{1j}$  as well as  on  $\sup_{t\in [0,T]} \|\xi_t^i\|_{TV}$ and  $\sup_{t\in
[0,T]} \|\tilde{\xi}_t^i\|_{TV}$. The latter suprema are in turn controlled  by  $e^{{\bar r}_i T}  \|\xi_0^i\|_{TV} $ by Remark \ref{bound}. By Gronwall's lemma we get $ \EE|X^i_{s,t}(x)-X^i_{s,t}(y)|^2\leq
C |x-y|^2$, which easily yields
\begin{equation*}
|P_{s,t}^i\varphi(x)-{P}^i_{s,t}\varphi(y)|^2\leq
C\|\varphi \|^2_{{\cal L}{\cal B}}|x-y|^2
\end{equation*}
as required.

{\it b)}  For notational simplicity we consider first the case $b=0$. Using similar types of inequalities as before, we have for all $s\leq t\leq T$, 
\begin{equation*}
\begin{split}
\EE|X^i_{s,t}(x)-\tilde{X}^i_{s,t}(x)|^2\leq  & ~
C'\bigg(\int_s^t \EE|X^i_{s,r}(x)-\tilde{X}^i_{s,r}(x)|^2 +  \EE| G^{i1}*\xi^1_s(X^i_{s,r}(x))-G^{i1}*\tilde{\xi}_s^1(\tilde{X}^i_{s,r}(x))|^2\\ 
&  \qquad + \EE|  G^{i2}* \xi^2_s(X^i_{s,r}(x))- G^{i2}*\tilde{\xi}^2_s(\tilde{X}^i_{s,r}(x))|^2 \,  dr \bigg) \\
\leq & C''\bigg(\int_s^{t} \EE|X^i_{s,r}(x)-\tilde{X}^i_{s,r}(x)|^2 + \EE \left|\int G^{i1}(X^i_{s,r}(x)-y)
(\xi_s^1(dy)-\tilde{\xi}^1_s(dy))\right|^2 \\
& \qquad + \EE\left|\int G^{i2}(X^i_{s,r}(x)-y)
(\xi_s^2(dy)-\tilde{\xi}^2_s(dy))\right|^2dr \bigg)\\
\end{split}
\end{equation*}
Since the functions $y\mapsto G^{ij}(X^i_s(x)-y)$ are  uniformly 
Lipschitz continuous, we deduce with
Gronwall's lemma that
\begin{equation*}\begin{split}
\EE(|X^i_{s,t}(x)-\tilde{X}^i_{s,t}(x)|^2)\leq & C \int_s^t
\|\xi^1_{r}-\tilde{\xi}^1_{r}\|^2_{{\cal L}{\cal B}^*}+\|\xi^2_r-\tilde{\xi}^2_r\|^2_{{\cal
L}{\cal B}^*}dr\\
\end{split}
\end{equation*}
which allows us to easily  conclude.  The case $b\not =0$ is similar with additional terms involving the kernels $H^{ij}$.

 \Box

\begin{lem}\label{PFK}
For each $t\geq 0$ and $\varphi\in C^3_0(\RR^d)$, the function $(s,x)\mapsto P^i_{s,t}\varphi (x)$ satisfies  the decay   condition \eqref{growthf} on each interval $[0,T']$ with $T'<t$. 
\end{lem}\label{PFK}

{\bf Proof:}   It is enough to show that 
\begin{equation}\label{controlderivgrow}
\sup_{s\in [0,T'], x\in \RR^d, |y|\leq 1}|x| |y|^{-1}  |P^i_{s,t}\varphi(x) - P^i_{s,t}\varphi(x+y)|  <\infty. 
\end{equation}
Since $\varphi$ and its derivatives are supported on $B(0,R)$ for some $R>0$, and given that $a^i $ and $b^i$ are globally Lipschitz under $(\mathbf{H})$, we get by  applying It\^o's formula on $\varphi$ that, for some constant $C_R>0$ and every $x,y\in \RR^d$ with  $|y|\leq 1$, it holds
\begin{equation*}
\begin{split}
|x|  |P^i_{s,t}\varphi(x) - P^i_{s,t}\varphi(x+ & y)|\leq   C_R |y|   {\bf 1}_{ |x|\leq R+1} \\
& + C_R  |x| \int_s^t \EE( |X^i_{s,r}(x)-X^i_{s,r}(x+y)| {\bf 1}_{  \min\{| X^i_{s,r}(x)|,  | X^i_{s,r}(x+y)| \leq R} )dr
\end{split}
\end{equation*}
By Cauchy-Schwartz inequality and the  proof of the previous Lemma we deduce that
\begin{equation*}
|x|  |P^i_{s,t}\varphi(x) - P^i_{s,t}\varphi(x+  y)|\leq   C_R |y|   {\bf 1}_{ |x|\leq R+1} + C_R  |x||y| \int_s^t \PP^{1/2}(  \min\{| X^i_{s,r}(x)|,  | X^i_{s,r}(x+y)|\} \leq R )dr
\end{equation*}
Since for each $z\in \RR^d$  the coefficients of  $ X^i_{s,r}(z)   $  grow linearly, 
$\sup_{r\in [s,T'] }\EE( |X^i_{s,r}(z) -z|^2)$ is finite so that, for $|z|>R$ one has
$\PP(  \min\{| X^i_{s,r}(z)| \leq R ) \leq C_{T'}/(|z|-R)^2$.  This and the previous estimate easily yield \eqref{controlderivgrow}. 
\Box
\medskip

{\bf Proof of Proposition \ref{uniqueness}:} Consider $t\in [0,T]$,  $\varphi\in C^3_0(\RR^d)$ and $\varepsilon< T-t$.   By the
Feynmann-Kac formula  \citeaffixed{KaSh}{e.g}, the function
$f^{(t+\varepsilon)}(s,x)=\EE(\varphi(X_{s,t+\varepsilon}^i(x)))=P^i_{s,t+\varepsilon}\varphi(x)$ is the
unique classic (bounded) solution of the linear parabolic problem
$$\partial_s f^{(t+\varepsilon)}(s,x)+ a^1_{k l}(\cdot,G^{i1}* \xi_s^1,G^{i2}* \xi_s^2) \partial_{x_k x_l}  f^{(t+\varepsilon)}(s,x)+
  b^1_{k }(\cdot,H^{i1}* \xi_s^1,H^{i2}* \xi_s^2) \partial_{x_k}f^{(t+\varepsilon)}(s,x)  =0$$with final condition at
time $s=t+\varepsilon$ equal to $\varphi(x)$.  By Lemma \ref{PFK},  the function $P^i_{s,t+\varepsilon}\varphi(x)$ moreover satisfies condition  \eqref{growthf} on $[0,t]$. Replacing it in the  first equation in \eqref{cross-diff} on $[0,t]$  and letting $\varepsilon\to 0$, we see that $\xi^1$ satisfies
\begin{equation*}
\langle  \xi_t^1,\varphi \rangle  =   \langle
\xi_0^1,P^1_{0,t}\varphi\rangle + \int_0^t \int
\big(r_1(x) -C^{11}*\xi_s^1(x)-C^{12}*\xi_s^2(x)\big)P^1_{s,t}\varphi(x)
 \xi_s^1(dx)ds,
\end{equation*}
and for all  $\varphi\in C^3_0(\RR^d)$ and then (by an elementary approximation argument) for all  $\varphi \in 
{\cal L}{\cal B}(\RR^d)$, which we chose such that $\|\varphi\|\leq 1$.
In a similar way we get
 \begin{equation*}
\langle \tilde{\xi}_t^1,\varphi \rangle  =   \langle
\xi_0^1,\tilde{P}^1_{0,t}\varphi\rangle + \int_0^t \int
\big(r_1(x)-C^{11}*\tilde{\xi}_s^1(x)-C^{12}*\tilde{\xi}_s^2(x)\big)\tilde{P}^1_{s,t}\varphi(x)
 \tilde{\xi}_s^1(dx)ds
\end{equation*}
for all  $\varphi \in 
{\cal L}{\cal B}(\RR^d)$ with $\|\varphi\|\leq 1$. Consequently, 
\begin{equation*}
\begin{split}
\langle   \xi_t^1 & -\tilde{\xi}_t^1,     \varphi \rangle^2 \\  \leq &
  \langle   \xi_0^1,  (P^1_{0,t}-\tilde{P}^1_{0,t})\varphi\rangle^2  + C  \int_0^t   \Bigg\{ \left[\int
(P^1_{s,t}-\tilde{P}^1_{s,t})\varphi(x) \xi_s^1(dx)\right]^2
+\left[\int \tilde{P}^1_{s,t}\varphi(x) (\xi_s^1(dx)-
\tilde{\xi}_s^1(dx))\right]^2 \\
& +\left[\int
C^{11}*(\xi_s^1-\tilde{\xi}_s^1)(x) P^1_{s,t}\varphi(x)
\xi_s^1(dx)\right]^2 +\left[\int C^{11}*\tilde{\xi}_s^1(x)
(P^1_{s,t}\varphi(x)-\tilde{P}^1_{s,t}\varphi(x))
\xi_s^1(dx)\right]^2  \\
& +\left[\int C^{11} *\tilde{\xi}_s^1(x)
\tilde{P}^1_{s,t}\varphi(x)
(\xi_s^1(dx)-\tilde{\xi}_s^1(dx))\right]^2+\left[\int
C^{12}*(\xi_s^2-\tilde{\xi}_s^2)(x) P^1_{s,t}\varphi(x)
\xi_s^1(dx)\right]^2  \\
&+\left[\int C^{12}*\tilde{\xi}_s^2(x)
(P^1_{s,t}\varphi(x)-\tilde{P}^1_{s,t}\varphi(x))
\xi_s^1(dx)\right]^2  +\left[\int C^{12}*\tilde{\xi}_s^2(x)
\tilde{P}^1_{s,t}\varphi(x)
(\xi_s^1(dx)-\tilde{\xi}_s^1(dx))\right]^2 \Bigg\} ds.
\\
\end{split}
\end{equation*}
Hence, 
\begin{equation*}
\begin{split}
\langle   \xi_t^1 & -\tilde{\xi}_t^1,     \varphi \rangle^2 \\
  \leq   &
C \sup_{y}|(P^1_{0,t}-\tilde{P}^1_{0,t})\varphi(y)|^2 + C  \int_0^t \bigg\{
\sup_{y}|(P^1_{s,t}-\tilde{P}^1_{s,t})\varphi(y)|^2 + \langle
\xi_s^1-\tilde{\xi}_s^1,\tilde{P}^1_{s,t} \varphi \rangle^2\\
&   + \sup_{y}| C^{11}*(\xi_s^1-\tilde{\xi}_s^1)(y)
|^2
   +\left[\int  C^{11}*\tilde{\xi}_s^1(x)
\tilde{P}^1_{s,t}\varphi(x)
(\xi_s^1(dx)-\tilde{\xi}_s^1(dx))\right]^2  \\ 
& +  \sup_{y}|  C^{12}*(\xi_s^2-\tilde{\xi}_s^2)(y)
|^2
   +\left[\int  C^{12}*\tilde{\xi}_s^2(x)
\tilde{P}^1_{s,t}\varphi(x)
(\xi_s^1(dx)-\tilde{\xi}_s^1(dx))\right]^2 
\bigg\}  ds, 
\end{split}
\end{equation*}
 for constants depending  on $\sup_{t\in [0,T]}
\|\xi_t^i\|_{TV}^2$, $\sup_{t\in [0,T]} \|\tilde{\xi}^i_t\|_{TV}^2$ and
$T$. The functions $x\mapsto C^{1j}(x-y)$ are Lipschitz
continuous uniformly in $y$ and, by Lemma
\ref{bornlip} {\it a)}, $\tilde{P}^1_{s,t} \varphi $ and   $C^{1j}*\tilde{\xi}_s^j(x)
\tilde{P}^1_{s,t}\varphi(x)$ are Lipschitz  continuous bounded
functions, uniformly in $s,t\in [0,T]$. Together with  Lemma
\ref{bornlip}  {\it b)},  this entails 
\begin{equation*}
\langle \xi_t^1-\tilde{\xi}_t^1,\varphi \rangle^2 \leq   ~  C
\int_0^t \sum_{i=1,2} \|\xi^i_s-\tilde{\xi}^i_s\|^2_{{\cal L}{\cal B}^*}
ds,
\end{equation*}
and we can analogously obtain a similar bound for $\langle \xi_t^2-\tilde{\xi}_t^2,\varphi \rangle^2 $.  Taking $\sup_{\|\varphi\|\leq 1}$, summing the two obtained inequalities and using Gronwall's lemma we
conclude that
$$\|\xi^i_t-\tilde{\xi}^i_t\|^2_{{\cal L}{\cal B}^*}=0$$
for all $t\in [0,T]$ and $i=1, 2$,  hence uniqueness for system \eqref{cross-diff}. \Box

\section{Regularity of the stochastic  flow and function solutions}

We next show  under two types of suitable assumptions on the coefficients and the initial condition that the solution $\xi^i_t(dx)$ has
a density $\xi^i_t(x)$ with respect to Lebesgue measure  for $i=1,\dots,M$. 

\begin{lem}\label{domination}
Let  $(\xi^1_t,\dots, \xi^M_t)_{t\in [0,T]}$ be the measure  solution of \eqref{cross-diff}. Then, for all $t\in [0,T]$,
$$ \xi_t^i   \leq  e^{\bar{r}_i T}  m_t^i,$$
where $m_t^i$ is the finite measure defined by 
\begin{equation}\label{m}
\langle m_t^i, \varphi \rangle:=  \EE\left(\int_{\RR^d}
\varphi(X_{0,t}^i(x))\xi_0^i(dx) \right)
\end{equation}
 for any bounded function $\varphi$.
\end{lem}

{\bf Proof:} We write the proof for $M=2$ and $i=1$, and omit  for notational simplicity the superscript $1$ in the flow $X^1_{s,t}(x)$.  Taking in the first equation  in \eqref{cross-diff}  the function $$f^{(t)}(s,x):=\EE\left(\varphi(X_{s,t}(x))\exp\left\{\int_s^t r_1(X_{s,r}(x))-C^{11}*\xi_r^1(X_{s,r}(x))-C^{12}*\xi_r^2(X_{s,r}(x)) dr\right\}\right),$$
which is by the Feynman-Kac formula \citeaffixed{KaSh}{see}   the unique classic (bounded) solution of the parabolic problem
\begin{equation*}
\begin{split}
0= & 
\partial_s f^{(t)}(s,x)+ \frac{1}{2} a^1_{k l}(\cdot,G^{i1}* \xi_s^1,G^{i2}* \xi_s^2) \partial_{x_k x_l}  f^{(t)}(s,x)+
  b^1_{k }(\cdot,H^{i1}* \xi_s^1,H^{i2}* \xi_s^2) \partial_{x_k}f^{(t)}(s,x)   \\
&   + \big(r_1(x)-C^{11}*\xi_s^1-C^{12}*\xi_s^2\big) f^{(t)}(s,x)\\
  \end{split}
  \end{equation*}
  with final condition at
time $s=t$ equal to $\varphi(x)$, we get that
 \begin{multline}
 \label{xiphieaT}
\langle   \xi_t^1  ,\varphi \rangle   =     \\  \EE\left( \int_{\RR^d} \varphi(X_{0,t}(x))\exp\left\{\int_0^t r_1(X_{0,s}(x)) -C^{11}*\xi_s^1(X_{0,s}(x))-C^{12}*\xi_s^2(X_{0,s}(x))ds\right\} \xi_0^1(dx) \right)  
\end{multline} 
for each  continuous bounded function
$\varphi\geq 0$. The above function $f^{(t)}(s,x)$ can be seen to satisfy condition \eqref{growthf} when $\varphi$ is as in Lemma \ref{PFK} and   $r_i$ is replaced by a bounded and Lipschitz continuous function $\tilde{r}_i$;  one then gets \ref{xiphieaT} up to an error that goes to $0$ 
 if $\tilde{r}_i$ goes to $r_i$ pointwise. It is the easy to deduce  \ref{xiphieaT} for $\varphi\in C_b(\RR^d)$. 

This yields  $\langle \xi_t^1,\varphi \rangle  \leq  e^{\bar{r}_1 T} \langle m_t^1,\varphi\rangle$ for all bounded continuous $\varphi$.  The measure  $\xi_t^1+e^{\bar{r}_1 T}m_t^1$ being regular, for each Borel set  $A$ and  $\epsilon >0$ we can find a closed set $B\subseteq A$ s.t. $\langle \xi^1_t - e^{\bar{r}_1 T} m^1_t,A\backslash B \rangle  \leq \epsilon$. Since the sequence of bounded continuous functions $f_k(x)=(1-kd(x,B))\vee 0$  pointwise converges  to $\mathbf{1}_B$ as $k\to \infty$, we have $\langle \xi^1_t-e^{\bar{r}_1 T} m^1_t, \mathbf{1}_B-f_k\rangle \to 0$ by dominated convergence w.r.t. the positive finite measure $|\xi_t^1- e^{\bar{r}_1 T} m_t^1 |$ .  It follows that for $k$ sufficiently large
\begin{equation*}
\langle \xi_t^1-e^{\bar{r}_1 T} m_t^1, A\rangle \leq  \langle \xi_t^1+m_t^1, A\backslash B  \rangle + \langle \xi_t^1-e^{\bar{r}_1 T} m_t^1,
   \mathbf{1}_B-f_k\rangle +  \langle \xi_t^1-e^{\bar{r}_1 T} m_t^1, f_k\rangle \leq 
   2\varepsilon ,
   \end{equation*}
  that is,  $\langle \xi_t^1,A \rangle  \leq  e^{\bar{r}_1 T} \langle m_t^1, A\rangle$.  
   \Box
   
   \bigskip
   
We immediately deduce from  \eqref{m} the following

\begin{coro} For any initial finite measure  $(\xi^1_0, \cdots, \xi^M_0)$ and 
in the uniform elliptic case: $\exists \lambda_i>0$ such that  $y^* a^i(x,v)y\geq \lambda_i| y |^2 \, \,  \forall x,y\in \RR^d, v\in \RR_+^M$,  the measure $\xi^i_t(dx)$ has a density $\xi^i_t(x)$ with respect to Lebesgue measure for all $t\in (0, T]$.
\end{coro}

Indeed, in that case, the law of the random variable $X_{0,t}^i(x)$ has a density with respect to Lebesgue measure. Lemma \ref{domination} allows us to conclude.

\me

As pointed out in Remark \ref{examples}, some natural  biological examples are not covered by this ellipticity assumption. We will next provide a finer result covering some non elliptic cases under additional regularity assumptions.  In  all the sequel, we assume
 
 \bigskip
 
 $(\mathbf{H})'$ :  Hypothesis  $(\mathbf{H})$ holds and moreover
 \begin{itemize}
 \item[i)]  $ \sigma^i(x,v_1,\dots,v_M)$ and $b^i(x,v_1,\dots,v_M)$  are respectively   ${\cal C}^{2,\alpha}(\RR^d\times [0,\infty)^M)$ and  ${\cal C}^{1,\alpha}(\RR^d\times [0,\infty)^M)$ for some $\alpha \in (0,1)$.
 \item[ii)] The  functions  $(G^{ij})_{1\leq i,j\leq M}$ and $( H^{ij})_{1\leq i,j\leq M}$  are   respectively  of class ${\cal C}^{2,\alpha}_b(\RR^d)$ and ${\cal C}^{1,\alpha}_b(\RR^d)$ for some $\alpha \in (0,1)$.
 \end{itemize}
 
 \medskip
 
 \begin{rem}\label{coefflowplusreg} Under assumptions $(\mathbf{H})'$ and by Remark \ref{bound},   $ \sigma(i,t,x)$ and $b(i,t,x)$ in \eqref{coefflow} are respectively  ${\cal C}^{2,\alpha}(\RR^d)$ and ${\cal C}^{1,\alpha}(\RR^d)$ for some $\alpha \in (0,1)$, uniformly in $[0,T]$.
\end{rem}

\begin{prop}\label{densityflow}
Assume   hypothesis  $(\mathbf{H})'$. If for some type $i$  the measure $\xi^i_0$ has a density, then  $\xi^i_t$ has a density for all $t\in [0, T]$.
\end{prop}

The proof will require classical regularity properties of  stochastic flows stated by \citeasnoun{Ku} and summarized in the next Lemma.

\begin{lem}\label{diffeoflow}
Under assumptions $(\mathbf{H})'$,  the process  $(s,t,x)\mapsto X^i_{s,t}(x)$ has a continuous version such that, a.s.  for each $s<t$ the function $x\mapsto X^i_{s,t}(x)$ is a global diffeomorphism of class ${\cal C}^{1,\beta}$ for all $\beta\in (0,\alpha)$.

Moreover, for each  $(t,y) \in [0,T]\times \RR^d$  the inverse mappings $\eta^i_{s,t}(y):= (X^i_{s,t})^{-1}(y)$, $0\leq s<t \leq T$, satisfy the stochastic differential equation

\begin{equation}\label{stochflow-1}
\eta^i_{s,t}(y)= y-\int_s^t \sigma(i,r,\eta^i_{r,t}(y)) \widehat{d}{B}_r^i - \int_s^t \hat{b}(i,r,\eta^i_{r,t}(y)) dr 
\end{equation}
where  for $k\in \{1,\cdots,d\}$, $\hat{b}_{k}(i,r,y)=b_{k}(i,r,y)- \sum_{q,l=1}^d \sigma_{lq}(i,r,y) \partial_{y_l} \sigma_{kq}(i,r,y)$ and $\widehat{d}{B}^i$ refers to the backward  It\^o integral with respect to the Brownian motion $B^i$  \citeaffixed{Ku}{see p. 194 in}.  

\medskip

Finally, for each  $(t,y) \in [0,T]\times \RR^d$  the (invertible)  Jacobian matrix $\nabla_y \eta^i_{s,t}(y)=\left( \partial_{y_l}\eta_{s,t}^{i,k}(y) \right)_{k,l=1}^d$ of $\eta^i_{s,t}(y)$ satisfies the system of backward linear stochastic differential equations :

\begin{equation}\label{derivflow}
 \partial_{y_l}\eta_{s,t}^{i,k}(y)= \delta_{kl} -\sum_{p,q=1}^d \int_s^t \partial_{x_p}      \sigma_{kq} (i,r,\eta^i_{r,t}(y)) \partial_{y_l}\eta_{r,t}^{i,p}(y) \widehat{d}{B}^{i,q}_r
 - \sum_{p=1}^d  \int_s^t \partial_{x_p}      \hat{b}_{k} (i,r,\eta^i_{r,t}(y)) \partial_{y_l}\eta_{r,t}^{i,p}(y) dr.
\end{equation}
\end{lem}

{\bf Proof:} Thanks to  assumption  $(\mathbf{H})'$ and Remark \ref{coefflowplusreg} we can apply Theorems 3.1 and 6.1 in Ch. II of  \citeasnoun{Ku}.

\Box

\bigskip

{\bf Proof of Proposition \ref{densityflow}:}  Let $\varphi\geq 0$ be a  bounded  measurable function in $\RR^d$. 
By the previous lemma, we can do the change of variable $X_{0,t}(x)=y$ in the integral inside the expectation in \eqref{m}, to get
$$  \langle m_t, \varphi \rangle =  \EE  \left(\int_{\RR^d} \varphi(y)  \xi_0^1(\eta_{0,t}(y)) | det \nabla \eta_{0,t}(y)| dy\right)  < +\infty .$$
Everything being positive in the above expression, Fubini's theorem yields that $m_t$ has the (integrable) density $y\mapsto \EE \left[ \xi_0^1(\eta_{0,t}(y)) | det \nabla \eta_{0,t}(y)|\right]   $ with respect to Lebesgue measure, and we conclude by Lemma \ref{domination}. 

\Box

\section{Convergence to local  competition} 

Our aim in this section is to describe some situations where the interaction range of the competition is much smaller than the one inducing  spatial dispersal. For example one may assume that animals interact for sharing resources as they are on the same place but diffuse  depending on the densities of the different species staying around  them in a larger neighbourhood. 
To model such situation, we suppose now that $C^{ij}=c^{ij} \gamma_{\varepsilon}$ for $c^{ij}\geq 0$ some fixed constant and  $\gamma_{\varepsilon}$ a suitable smooth approximation of the Dirac mass as $\varepsilon\to 0$. Our goal is to show that, under additional regularity assumptions, the (unique) solution $\xi=(\xi^{1}, \cdots, \xi^{M})$ of Equation \eqref{cross-diff} given by Theorem  \ref{existunique} for such competition coefficients converges, as $\varepsilon\to 0$, to a weak function solution of the system of Equations \eqref{cross-diffusion2}.

\medskip

In what follows,   stronger conditions on the coefficients will be enforced, namely:

\medskip

 $(\mathbf{H})''$ :  Hypothesis  $(\mathbf{H})$ holds and moreover
 \begin{itemize}
 \item[i)]  $ \sigma^i(x,v_1,\dots,v_M)$ and $b^i(x,v_1,\dots,v_M)$  are respectively   ${\cal C}^{3,\alpha}(\RR^d\times [0,\infty)^M)$ and  ${\cal C}^{2,\alpha}(\RR^d\times [0,\infty)^M)$ for some $\alpha \in (0,1)$. 
 \item[ii)] Functions  $(G^{ij})_{1\leq i,j\leq M}$ and $( H^{ij})_{1\leq i,j\leq M}$  are respectively ${\cal C}^{3,\alpha}(\RR^d)$ and ${\cal C}^{2,\alpha}(\RR^d)$ for some $\alpha \in (0,1)$. Moreover, functions  $(C^{ij})_{1\leq i,j\leq M}$ are integrable in $\RR^d$ and  have bounded derivatives.
 \item [iii)] Functions $r_i$  have bounded derivatives.
 \end{itemize}

 \begin{rem}\label{coefflowveryreg} Under assumptions $(\mathbf{H})''$,  functions $ \sigma(i,t,x)$ and $b(i,t,x)$ in \eqref{coefflow} are respectively  $C^{3,\alpha}_b(\RR^d)$ and $C^{2,\alpha}(\RR^d)$ for some $\alpha \in (0,1)$, uniformly in $[0,T]$, and with bounds that do not depend on the  kernels $C^{ij}$  (cf. Remark \ref{bound}).
\end{rem}

We will next establish

\begin{theo}\label{existlocalcompet}
Assume  that  hypothesis $(\mathbf{H})''$ hold,  and that the initial measures $(\xi^1_0, \cdots, \xi^M_0)$ have  densities in $L^1\cap L^{\infty}$ and   distributional derivatives  in $L^{\infty}$. 

\medskip

Furthermore, assume that $C^{ij}=c^{ij} \gamma_{\varepsilon}$ for $c^{ij}\geq 0$ some fixed constant and  $\gamma_{\varepsilon}=\gamma(x/\varepsilon) \varepsilon^{-d}$ for  some regular function  $\gamma\geq 0$ satisfying $\int_{\RR^d}\gamma(x)dx=1$ and $\int_{\RR^d}|x|\gamma(x)dx<\infty$. 

\medskip

Then, for each $T>0$ the unique weak function-solution  $\xi^{\varepsilon}$ to Equation \eqref{cross-diff} converges in the space $C([0,T],{\cal M}^M)$ (endowed with the uniform topology) at speed $\varepsilon$ with respect to  the dual Lipschitz norm, to a solution $u=(u^1,\cdots, u^M)$ of the   non local cross-diffusion system with local competition: 
\begin{eqnarray}
\label{cross-diff-localcompet}
\langle u_t^i,f^{i}_t \rangle & =  & \langle
\xi_0^i,f^{i}_0\rangle  +  \int_0^t \int \bigg\{\frac{1}{2} \sum_{k,l}a^i_{k l}(\cdot, G^{i1}* u_t^1,\cdots, G^{iM}* u_t^M)\, \partial_{x_k x_l}  f^{i}_s  
\notag
\\ &&  \hspace{0.1cm}+
\sum_{k}  b^i_{k }(\cdot,H^{i1}* u_t^1,\cdots, H^{iM}* u_t^M)\, \partial_{x_k}   f^{i}_s +\big(r_{i}- \sum_{j=1}^M c^{ij}u_s^j\big)f^{i}_s+\partial_s f^i_s
\bigg\} (x)\ u_s^i(x)dxds.\notag\\
\ 
\end{eqnarray}
Moreover the function $u$ is the unique  function solution of \eqref{cross-diff-localcompet} such that 
\be
\label{cond-uni}
\sup_{t\in [0,T]} \| u_t\|_1 +  \| u_t\|_{\infty} +   \| \nabla  u_t\|_{\infty}   <+\infty.\ee
\end{theo}

To prove Theorem \ref{existlocalcompet}, we will extend   to a  convergence argument some of the techniques previously used in the uniqueness result. The same dual norm will be used, along with  some additional estimates and technical results:

\begin{lem}\label{estimagrad}
Under the assumptions of Theorem \ref{existlocalcompet},  for each $t\in [0, T]$  the functions $(\xi^\varepsilon_t)^i, i=1,\dots, M$ have bounded first order derivatives. Moreover, there exists for each $i$ a constant $K_i>0$ depending on the functions $C^{ij},j=1,\cdots,M$ only through their $L^1$ norms $c^{ij}$ (and not on $\varepsilon$), such that
\begin{equation*}
\max\{ \sup_{t\in [0,T]} \| (\xi^\varepsilon_t)^i\|_{\infty},
\sup_{t\in [0,T]} \| \nabla(\xi^\varepsilon_t)^i\|_{\infty} \}<K_i ,\quad \forall\, i=1,\dots M.
\end{equation*}
\end{lem}
This result relies on  an enhancement of Lemma \ref{diffeoflow} needing finer properties of stochastic flows established by  \citeasnoun{Ku1}.

\begin{lem}\label{Jacflow}
Under  assumption $(\mathbf{H})''$ 
  for each $i=1,\dots,M$ and $p\geq 2$,  there exist finite constants $K_{i1}(p)>0$ and $K_{i2}(p)>0$ not depending on the  kernels $C^{ij}$ such that  for all $t\in [0, T]$,
$$\sup_{y \in \RR^d}\EE\left( \sup_{s\in [0,t]} |  \nabla_y \eta^i_{s,t}(y) |^p\right)<K_{i1}(p) \mbox{  and } \quad 
\sup_{y \in \RR^d}\EE \left( \sup_{s\in [0,t]} | det  \nabla_y \eta^i_{s,t}(y)|^p  \right) <K_{i2}(p). $$
Moreover, for each $s<t$ with  $s,t\in [0,T]$ the function $y\mapsto det  \nabla_y \eta^i_{s,t}(y)$ is a.s. differentiable and there exists $K_{i3}(p)>0$ not depending on the  kernels $C^{ij}$ such that  for all $t\in [0, T]$,
$$\sup_{y \in \RR^d}\EE \left( \sup_{s\in [0,t]} |\nabla_y  \left[ det  \nabla_y \eta^i_{s,t}(y)\right]|^p  \right) <K_{i3}(p). $$
\end{lem}  

{\bf Proof of Lemma \ref{Jacflow}:}
For fixed $i\in \{1,\dots,M\}$  and $t\in [0,T]$ we define  coefficients $\beta:[0,t]\times \RR^d\to \RR^d$ and $ A^q:[0,t]\times \RR^d\to \RR^d$, $q=1,\dots, d$ with components $\beta_k $ and $A^q_k$, $k=1,\dots,d $ by 
$$\beta_{k}(s,y):= -\hat{b}_{k}(i,t-s,y), \quad A^q_{k}(s,y):=\sigma_{k q}(i,t-s,y)$$
(see Lemma \ref{diffeoflow} for the notation)
and the process $(Z_s(y);  s\in [0,t],y \in\RR^d)$ by $Z_s(y)=\eta_{t-s,t}^i(y)$. Then, denoting by $W=(W^1,\dots,W^d)$ the standard $d-$ dimensional Brownian motion $W_s:=B^i_{t-s}-B^i_t$, it easily follows from Lemma \ref{diffeoflow} that $Z_s(y)$ satisfies the classic It\^o stochastic differential equation 
\begin{equation}\label{stochflow-reverted}
Z_{s}(y)= y+\int_0^s A(r, Z_r(y) )d W_r +  \int_0^s \beta(r, Z_r(y)) dr \, ,
\end{equation}
whereas the associated Jacobian matrix satisfies the linear system:
\begin{equation}\label{derivflowrevert}
 \partial_{y_l}Z_s^k(y)= I_d +\sum_{m,q=1}^d \int_0^s \partial_{x_p}      A_{k}^q (r,Z_r(y)) \partial_{y_l}Z_r^{m}(y) d
 W^{q}_r
 + \sum_{m=1}^d  \int_0^s \partial_{x_p}      \beta_{k} (r,Z_r(y)) \partial_{y_l}Z_r^{m}(y) dr.
\end{equation}
Notice that the (non-homogenous) coefficients of this linear SDE are  uniformly bounded (cf. Remark \ref{coefflowveryreg}) independently of  kernels $C^{ij}$.
 Using the  Burkholder-Davis-Gundy inequality, the boundedness of the  derivatives of $A^q$ and $\beta$ and Gronwall's lemma, we deduce that
 \begin{equation}\label{gradZ}
 \EE\left( \sup_{s\in [0,t]} |  \nabla_y Z_s(y) |^p\right)<K_{i1}(p)
 \end{equation}
for some constant $K_{i1}(p)$ which depends on bounds for those derivatives and on $\sup_{t\in [0,T]} \|\xi_t^i\|_{TV}$ (cf. Remark \ref{coefflowplusreg}) but does not depend on $y\in \RR^d$.  This yields the first asserted estimate. 

In order to get the estimates for the determinant and its gradient, we rewrite  \eqref{stochflow-reverted} in Stratonovich form\begin{equation*}
Z_{s}(y)= y+\int_0^s A(r, Z_r(y) )\circ d W_r +  \int_0^s \beta^{\circ}(r, Z_r(y)) dr 
\end{equation*}
where $\beta^{\circ}(r,x)=\beta(r,x)-\frac{1}{2} \sum_{l,q} A^q_l(r,x)\partial_{x_l}A^q(r,x)$. By the proof of Lemma 4.3.1 of  \citeasnoun{Ku1}, $det \nabla_y Z_{s}(y)$ satisfies the linear Stratonovich stochastic differential equation 
\begin{equation}\label{eqdet}
det \nabla_y Z_{s}(y)=1+ \int_0^s  det \nabla_y Z_r(y) \sum_{k=1}^d\left[ \sum_{q=1}^d \partial_{y_k}A_k^q(r,Z_r(y))\circ d W^q_r +
\partial_{y_k}\beta^{\circ}_k(r,Z_r(y)) dr\right].
\end{equation}
Again,  the  coefficients of this scalar linear SDE are  uniformly bounded  independently of the  kernels $C^{ij}$. 
Using Burkholder-Davis-Gundy inequality in the It\^o's form of the previous equation, we deduce using also Gronwall's lemma that 
\begin{equation}\label{detgradZ}
\EE\left( \sup_{s\in [0,t]} |  det \nabla_y Z_s(y) |^p\right)<K_{i2}(p)
\end{equation}
for some constant $K_{i2}(p)$ depending on bounds on the (up to second order) derivatives of $\sigma^i$ and  (up to first order derivatives) of $b$, on $\sup_{t\in [0,T]} \|\xi_t^i\|_{TV}$ and on the constant $C_M$ in assumption $(\mathbf{H})''$ i). This  yields the second required estimate. 
Remark \ref{bound} ensures that the constants $K_{i1}(p)$ and $K_{i2}(p)$ do not depend on the kernels $C^{ij}$ nor on $\varepsilon$. 

Finally, under  assumptions $(\mathbf{H})''$  we deduce   from equation \eqref{eqdet} and Theorem 3.3.3 of  \citeasnoun{Ku1} (see also Exercise 3.1.5 therein)   the a.s. differentiability of  the mapping $y\mapsto det  \nabla_y \eta^i_{s,t}(y)$, and the fact  that 
 its derivative with respect to $y_l$ satisfies
\begin{equation*}
\begin{split}
\partial_{y_l}&[ det \nabla_y Z_{s}(y)]=  \int_0^s  \partial_{y_l} [det \nabla_y Z_r(y)] \sum_{k=1}^d\left[ \sum_{q=1}^d \partial_{y_k}A_k^q(r,Z_r(y))\circ dW^q_r +
\partial_{y_k}\beta^{\circ}_k(r,Z_r(y)) dr \right] \\
&+ \int_0^s  det \nabla_y Z_r(y) \sum_{m,k=1}^d \partial_{y_l}Z^m_s(y)\left[ \sum_{q=1}^d \partial^2_{y_m y_k}A_k^q(r,Z_r(y))\circ d W^q_r +
\partial_{y_m y_k}\beta^{\circ}_k(r,Z_r(y))dr \right] .
\end{split}
\end{equation*}
Note that all coefficients inside the square brackets are uniformly bounded functions (independently of  kernels $C^{ij}$). 
Writing this equation in It\^o's form, we now deduce with the Burkholder-Davis-Gundy inequality that $\phi(s):=\EE \left( \sup_{r\in [0,s]} |\nabla_y  \left[ det  \nabla_y Z_r(y)\right]|^p  \right)$ satisfies the inequality
$$
\phi(s)\leq C' \int_0^s  \phi(r)dr + C'' \int_0^s \EE \left( \sup_{\theta \in [0,r]}  | det  \nabla_y Z_{\theta}(y)|^p   \sup_{\theta\in [0,r]} |  \nabla_yZ_{\theta}(y) |^p  \right)dr.$$
By Cauchy-Schwarz inequality and the estimates \eqref{gradZ} and  \eqref{detgradZ} with $2p$ instead of $p$, the above expectation is seen to be bounded uniformly in $r\in [0, t], y\in \RR^d$. We deduce by Gronwall's lemma that 
\begin{equation*}
\EE \left( \sup_{s\in [0,t]} |\nabla_y  \left[ det  \nabla_y Z_s(y)\right]|^p  \right)<K_{i3}(p)
\end{equation*}
for some constant $K_{i3}(p)$ as required, and conclude the third asserted estimate.

\Box

\bigskip

{\bf Proof of Lemma \ref{estimagrad}:} We again consider $M=2$, $i=1$ and omit  the superscript $1$  in the process  $X^1_{s,t}(x)$, the inverse flow and its derivative.   By  Lemma \ref{diffeoflow}  we can do the change of  variables $X_{0,t}(x)=y$ in the integral with respect to $dx$ inside the expectation in  \eqref{xiphieaT}. Using the semigroup property of the flow and its inverse  stated by \citename{Ku}  \citeyear{Ku,Ku1}  together with  Fubini's theorem (thanks to Lemma {\ref{Jacflow}),  we deduce  that for a.e. $y\in \RR^d$, 
\begin{equation}\label{xitexplicit}
\xi_t^1(y)= \EE\left[ \Psi(t,y)\right]
\end{equation}
where $ \Psi(t,y)$ is the random function
\begin{equation*}
 \Psi(t,y) : =  \exp\left\{\int_0^t(r_1(\eta_{r,t}(y))-C^{11}*\xi_r^1(\eta_{r,t}(y))-C^{12}*\xi_r^2(\eta_{r,t}(y)))dr\right\} \xi_0^1(\eta_{0,t}(y)) det \nabla_y\eta_{0,t}(y).
  \end{equation*}
  Notice that we have used the fact that $det \nabla_y\eta_{0,t}(y)>0$, which follows from $det \nabla_y\eta_{r,t}(y)\not = 0$  for all $r\in [0,t]$ and    $r\mapsto \nabla_y\eta_{r,t}(y)$ being continuous with  value  $I_d$ at $r=t$.  The bound on $\sup_{t\in [0,T]} \| \xi_t^i\|_{\infty}$  readily follows from the previous identity, the assumptions on $\xi_0^i$ and the second estimate in Lemma \ref{Jacflow}.
  
  The function $ y\mapsto \Psi(t,y)$ is moreover continuously differentiable,  by Lemmas \ref{diffeoflow} and \ref{Jacflow}. 
   Since the kernels $C^{11}$ and $C^{12}$ have bounded derivatives we deduce that,  a.s.  
\begin{equation}\label{gradxi1}
\begin{split}
\nabla \Psi(t,y)  = &   \exp\left\{\int_0^t r_1(\eta_{r,t}(y))-C^{11}*\xi_r^1(\eta_{r,t}(y))-C^{12}*\xi_r^2(\eta_{r,t}(y)) dr\right\}   \\
& \,\times \bigg[
 \xi_0^1(\eta_{0,t}(y)) det \nabla \eta_{0,t}(y)  \int_0^t    \nabla^* \eta_{r,t}(y)   \left[  \nabla r_1 - (\nabla C^{11}) *\xi_r^1-(\nabla C^{12})*\xi_r^2  \right] (\eta_{r,t}(y))dr  \\ 
& \qquad +  \nabla^* \eta_{0,t}(y)\nabla\left[ \xi_0^1\right] (\eta_{0,t}(y)) det \nabla \eta_{0,t}(y) + \xi^1_0( \eta_{0,t}(y))  \nabla  \left[ det  \nabla \eta_{0,t}(y)\right] \bigg] \\
\end{split}
\end{equation}
for all $y \in \RR^d$. From  Lemma \ref{Jacflow}, thanks to Cauchy-Schwarz inequality we get  that 
$$\EE \left( |  det \nabla \eta_{0,t}(y) | \int_0^t | \nabla \eta_{r,t}(y)|dr + | \nabla \eta_{0,t}(y)| \, |  det \nabla \eta_{0,t}(y)| + | \nabla  \left[ det  \nabla \eta_{0,t}(y) \right]  | \right)<\infty.$$
Thus,  we can take derivatives inside the expectation  \eqref{xitexplicit}
 and  deduce the existence of 
$$\nabla \xi_t^1(y)= \EE\left[ \nabla \Psi(t,y)\right],$$ and moreover that 
$\sup_{t\in [0, T]} \| \nabla  \xi_t^1 \|_{\infty}<\infty$.  Similarly,  $\sup_{t\in [0, T]} \| \nabla  \xi_t^2 \|_{\infty}<\infty$.  We can now rewrite \eqref{gradxi1} as
\begin{equation}\label{gradxi2}
\begin{split}
\nabla \Psi(t,y)  = &   \exp\left\{\int_0^t  r_1(\eta_{r,t}(y)) -C^{11}*\xi_r^1(\eta_{r,t}(y))-C^{12}*\xi_r^2(\eta_{r,t}(y))dr\right\}   \\
& \,\times \bigg[
\xi_0^1(\eta_{0,t}(y)) det \nabla \eta_{0,t}(y)  \int_0^t   \nabla^* \eta_{r,t}(y)   \left[  \nabla r_1 - C^{11}* \nabla  \xi_r^1- C^{12}*\nabla \xi_r^2  \right] (\eta_{r,t}(y))   dr  \\ 
& \qquad + \nabla^* \eta_{0,t}(y)\nabla\left[ \xi_0^1\right] (\eta_{0,t}(y)) det \nabla \eta_{0,t}(y) + \xi^1_0( \eta_{0,t}(y))  \nabla  \left[ det  \nabla \eta_{0,t}(y)\right]  \bigg]. \\
\end{split}
\end{equation}
Since $\| C^{1j}\|_1 =c^{1j}$, we have $  \| C^{1j}*\nabla \xi_s^j\|_{\infty}\leq c^{1j}  \| \nabla \xi_s^j\|_{\infty} $ for all $s\in [0, T]$. Taking expectation in \eqref{gradxi2} and using the estimates in Lemma \ref{Jacflow}, we deduce that for all $t\in [0,T]$, 
\begin{equation*}
\| \nabla \xi_t^1\|_{\infty} \leq C'' \int_0^t ( \| \nabla \xi_r^1\|_{\infty} +   \| \nabla \xi_r^2\|_{\infty}) dr+ C'''
\end{equation*}
for constants $C''',C''>0$  depending on the functions  $C^{1j}$  only through their $L^1$ norms $c^{1j}$ (in particular not depending on $\varepsilon$).  Summing the later estimate with the analogous one for $\| \nabla \xi_t^2\|_{\infty} $, we conclude thanks to Gronwall's lemma. 

\Box

\bigskip

We are now ready for the
\medskip

{\bf Proof of Theorem  \ref{existlocalcompet}:} Again, we write  the proof in the case $M=2$. Let $\varepsilon>\bar{\varepsilon}>0$. To lighten notation, we denote simply by $\xi=(\xi^1,\xi^2)$ and
$\bar{\xi}=(\bar{\xi^1},\bar{\xi^2})$  two solutions of
system \eqref{cross-diff} in $[0,T]$ respectively with $C^{ij}=c^{ij} \varphi_{\varepsilon}$ and $\bar{C}^{ij}:=c^{ij} \varphi_{\bar{\varepsilon}}$.
Proceeding as in the proof of Proposition \ref{uniqueness}, we deduce that for all  function $\varphi$ with $\| \varphi\|_{{\cal L}{\cal B}}\leq 1$,
\begin{equation}\label{boundsxibarxi}
\begin{split}
\langle   \xi_t^1 & -\bar{\xi}_t^1,     \varphi \rangle^2 \\  \leq &
  \langle   \xi_0^1,  (P^1_{0,t}-\bar{P}^1_{0,t})\varphi\rangle^2  + C  \int_0^t   \Bigg\{ \left[\int
(P^1_{s,t}-\bar{P}^1_{s,t})\varphi(x) \xi_s^1(x)dx\right]^2
+\left[\int \bar{P}^1_{s,t}\varphi(x) (\xi_s^1(x)-
\bar{\xi}_s^1(x))dx\right]^2 \\
& +\left[\int
C^{11}*(\xi_s^1-\bar{\xi}_s^1)(x) P^1_{s,t}\varphi(x)
\xi_s^1(x)dx\right]^2 +\left[\int C^{11}*\bar{\xi}_s^1(x)
(P^1_{s,t}\varphi(x)-\bar{P}^1_{s,t}\varphi(x))
\xi_s^1(x)dx\right]^2  \\
& +\left[\int C^{11} *\bar{\xi}_s^1(x)
\bar{P}^1_{s,t}\varphi(x)
(\xi_s^1(x)-\bar{\xi}_s^1(x))dx\right]^2\\
& +  \left[\int [ C^{11}- \bar{C}^{11}] *\bar{\xi}_s^1(x)
\bar{P}^1_{s,t}\varphi(x)
\bar{\xi}_s^1(x)dx\right]^2\\ 
& +\left[\int
C^{12}*(\xi_s^2-\bar{\xi}_s^2)(x) P^1_{s,t}\varphi(x)
\xi_s^1(x)dx\right]^2  +\left[\int C^{12}*\bar{\xi}_s^2(x)
(P^1_{s,t}\varphi(x)-\bar{P}^1_{s,t}\varphi(x))
\xi_s^1(x)dx\right]^2 \\
& +\left[\int C^{12}*\bar{\xi}_s^2(x)
\bar{P}^1_{s,t}\varphi(x)
(\xi_s^1(x)-\bar{\xi}_s^1(x))dx\right]^2  
\\
&+  \left[\int [ C^{12}- \bar{C}^{12}] *\bar{\xi}_s^2(x)
\bar{P}^1_{s,t}\varphi(x)
\bar{\xi}_s^1(x)dx\right]^2 \Bigg\} ds.\\
\end{split}
\end{equation}
Thanks to Lemma \ref{estimagrad}, $P^1_{s,t}\varphi\xi_s^1$ is a bounded Lipschitz function with  Lipschitz norm  bounded  independently of $\varepsilon$, $\bar{\varepsilon}$ and $s,t\in [0,T]$. We thus can rewrite and bound the first term in  the third  line of   of \eqref{boundsxibarxi} 
as follows:  $$\left[ \int
(\xi_s^1-\bar{\xi}_s^1)(y)  C^{11} *(P^1_{s,t}\varphi\xi_s^1) (y) \, dy\right]^2 \leq  C\|\xi^1_{s}-\bar{\xi}^1_{s}\|^2_{{\cal L}{\cal B}^*} $$
for some $C>0$ not depending on  $\varepsilon,\bar{\varepsilon}$.  The  second term in the third line   is controlled by
\begin{equation*}
  C\| \bar{\xi}_s^1\|^2_{\infty} \| \xi_s^1\|^2_{TV}  \sup_{x\in \RR^d}\left| 
P^1_{s,t}\varphi(x)-\bar{P}^1_{s,t}\varphi(x)\right|^2
\leq  C
 \| \bar{\xi}_s^1\|_{\infty}^2 \| \xi_s^1\|_{TV}^2  \int_s^t \|\xi^1_r-\bar{\xi}^1_r\|^2_{{\cal
L}{\cal B}^*}+ \|\xi^2_r-\bar{\xi}^2_r\|^2_{{\cal
L}{\cal B}^*}dr
\end{equation*}
 thanks to Lemma \ref{bornlip} b). The term in the fourth line is easily controlled by $C\|\xi^1_{s}-\bar{\xi}^1_{s}\|^2_{{\cal L}{\cal B}^*}$. Using  the fact that,  by Lemma  \ref{estimagrad},  $\bar{\xi}^1_s$ has  derivatives uniformly bounded independently of $\bar{\varepsilon}>0$ and $s\in [0,T]$,
 we deduce by the assumption on $\gamma$ that 
 $$\sup_{x\in \RR^d}|[C^{11}- \bar{C}^{11}] *\bar{\xi}_s^1(x)|^2\leq C |\varepsilon -\bar{\varepsilon}|^2. $$
 A similar upper  bound then follows for the term in the fifth line of \eqref{boundsxibarxi}. The last three lines can be bounded in a similar way, and the first line on the right hand side is bounded in terms of dual Lipschitz distances by similar arguments as in Proposition \ref{uniqueness}.  Proceeding  in a similar way as therein, we now obtain the estimate
 \begin{equation*}
\langle \xi_t^1-\bar{\xi}_t^1,\varphi \rangle^2 \leq   ~  C
\int_0^t \sum_{i=1,2} \|\xi^i_s-\bar{\xi}^i_s\|^2_{{\cal L}{\cal B}^*}
ds+ C|\varepsilon -\bar{\varepsilon}|^2,
\end{equation*}
where the constants do not depend on $t\in [0,T]$ nor on $\varepsilon$ or $\bar{\varepsilon}$. Taking suitable suprema, the latter estimate can thus be  strengthened  to
 \begin{equation*}
 \sup_{r\in [0,t]} \|\xi^1_r-\bar{\xi}^1_r\|^2_{{\cal L}{\cal B}^*} \leq   ~  C
\int_0^t \sum_{i=1,2}\sup_{r\in [0,s]} \|\xi^i_r-\bar{\xi}^i_r\|^2_{{\cal L}{\cal B}^*}
ds+ C|\varepsilon -\bar{\varepsilon}|^2,
\end{equation*}
which when summed with the corresponding estimate for $i=2$ yields 
 \begin{equation}\label{xibarxi}
 \sum_{i=1,2} \sup_{r\in [0,T]} \|\xi^i_r-\bar{\xi}^i_r\|^2_{{\cal L}{\cal B}^*} \leq C|\varepsilon -\bar{\varepsilon}|^2, 
\end{equation}
 after applying  Gronwall's lemma.  Therefore, as $\varepsilon$ goes to $0$, the sequence $(\xi^{\varepsilon})_{\varepsilon>0}$ is Cauchy  in the Polish space $C([0,T],{\cal M}^M)$ and thus converges to some element in that space. 
 Dunford-Pettis criterion for weak compactness in $L^1$, together with  the uniform bounds both in $L^1$ and $L^{\infty}$  for  $(\xi^{\varepsilon}_t)_{\varepsilon>0}$,  imply that the components of  the previous limit  have densities for each $t\in [0,T]$,  which we denote by $u^i_t(x)$, and which  satisfy the same $L^1$ and $L^{\infty}$ bounds. 

Weak $L^1$-convergence  is however not enough  to identify  $u$ as a solution of  \eqref{cross-diff-localcompet} and some regularity of the limit  will be needed to do so. 
 We denote by $\hat{P}_{s,t}^i(x,dy)$ the semigroup associated with the SDE with coefficients defined in terms of the measures  $(u_t^i(x)dx)_{t\in [0,T]}$ as previously. For  $\varphi$ such that  $\| \varphi\|_{{\cal L}{\cal B}}\leq 1$,  we set
$$\Psi^1(t,\varphi):= \langle u_t^1,\varphi \rangle  -   \langle
\xi_0^1,\hat{P}^1_{0,t}\varphi\rangle - \int_0^t \int
\big(r_1(x)-c^{11}u_s^1(x)-c^{12}u_s^2(x)\big)\hat{P}^1_{s,t}\varphi(x)
u_s^1(x)dx ds.$$
Since  $\xi=\xi^{\varepsilon}$ satisfies $$
\langle  \xi_t^1,\varphi \rangle -   \langle
\xi_0^1,P^1_{0,t}\varphi\rangle - \int_0^t \int
\big(r_1(x) -C^{11}*\xi_s^1(x)-C^{12}*\xi_s^2(x)\big)P^1_{s,t}\varphi(x)
 \xi_s^1(x)dxds=0,
$$
proceeding in a similar way as done to obtain the estimate \eqref{boundsxibarxi}, we deduce that
 \begin{multline}\label{sumofterms}
| \Psi^1(t,\varphi)|^2\leq  \sum_{i=1,2} \sup_{r\in [0,t]} \|\xi^i_r-u^i_r\|^2_{{\cal L}{\cal B}^*}  + \int_0^t\left[\int [ C^{1i}*u_s^i(x) -c^{1i}u_s^i(x)]
\hat{P}^1_{s,t}\varphi(x)
u_s^1(x)dx\right]^2 ds
\end{multline}
(the last term  corresponding to the sum of the fifth and last lines in   \eqref{boundsxibarxi} when $\bar{\varepsilon}=0$).   
Now,  by  Lemma \ref{estimagrad}, there exists a constant $K>0$ (independent of $\varepsilon>0$ and $s\in [0,T]$) such that for any $\varepsilon>0$ one has $ |\langle ( \xi^{\varepsilon}_{s})^i, \partial_{x_{l}} \varphi \rangle |\leq K \|\varphi\|_{1}$
for any   ${\cal C}^\infty$  compactly supported function $\varphi$. By letting $\varepsilon\to 0$,  the same bound is satisfied by $u^i_s$. By standard results on Sobolev spaces  \citeaffixed{Bre}{see e.g. Proposition IX.3 of} we get  that $u^i_s$ has distributional derivatives in $L^{\infty}$ and that $u^i_s$ is  Lipschitz continuous with  Lipschitz constant less than or equal to $K$.  Since $C^{1i}*u_s^i(x) -c^{1i}u_s^i(x)=c^{1i}\int \gamma(z)[ u^i_s(x+\varepsilon z)- u^i_s(x)]dz$
and $\int \gamma(z)|z|dz<\infty$, we deduce that $\| C^{1i}*u_s^i -c^{1i}u_s^i\|_{\infty} \leq C\varepsilon$, which combined with the bound
 $$\sum_{i=1,2} \sup_{r\in [0,t]} \|\xi^i_r-u^i_r\|^2_{{\cal L}{\cal B}^*}\leq C\varepsilon^2 \, ,$$
 following from \eqref{xibarxi},  yields  $| \Psi(t,\varphi)|\leq C\varepsilon $ for all $\varepsilon>0$. That is, $u$ solves  \eqref{cross-diff-localcompet}.

  \medskip
  
  Let us finally prove the uniqueness of the solution $u$. Recall that $\hat{P}_{s,t}^i(x,dy)$ denotes the associated diffusion semigroup. Consider a second function solution $v$  in $C([0,T],{\cal M}^M)$ satisfying \eqref{cond-uni} and with  associated semigroups denoted  $\check{P}_{s,t}^i(x,dy)$. Then, 
  \begin{equation}\label{bounduv}
\begin{split}
\langle   u_t^1 & -v_t^1,     \varphi \rangle^2 \\  \leq &
  \langle   u_0^1,  (\hat{P}^1_{0,t}-\check{P}^1_{0,t})\varphi\rangle^2  + C  \int_0^t   \Bigg\{ \left[\int
(\hat{P}^1_{s,t}-\check{P}^1_{s,t})\varphi(x) u_s^1(x)dx\right]^2
+\left[\int \check{P}^1_{s,t}\varphi(x) (u_s^1(x)-
v_s^1(x))dx\right]^2 \\
& +\left[\int
c^{11}(u_s^1-v_s^1)(x) \hat{P}^1_{s,t}\varphi(x)
u_s^1(x)dx\right]^2 +\left[\int c^{11}v_s^1(x)
(\hat{P}^1_{s,t}\varphi(x)-\check{P}^1_{s,t}\varphi(x))
u_s^1(x)dx\right]^2  \\
& +\left[\int c^{11} v_s^1(x)
\check{P}^1_{s,t}\varphi(x)
(u_s^1(x)-v_s^1(x))dx\right]^2\\
& +\left[\int
c^{12}(u_s^2-v_s^2)(x) \hat{P}^1_{s,t}\varphi(x)
u_s^1(x)dx\right]^2  +\left[\int c^{12}v_s^2(x)
(\hat{P}^1_{s,t}\varphi(x)-\check{P}^1_{s,t}\varphi(x))
u_s^1(x)dx\right]^2 \\
& +\left[\int c^{12}v_s^2(x)
\check{P}^1_{s,t}\varphi(x)
(u_s^1(x)-v_s^1(x))dx\right]^2   \Bigg\} ds.\\
\end{split}
\end{equation}
The functions  $\hat{P}^1_{s,t}\varphi u_s^1$ and $
\check{P}^1_{s,t}\varphi v_s^1$ having   Lipschitz norm  bounded   independently of  $s,t\in [0,T]$, the first term in  the third line and the term in the forth line above are bounded by $ C\|u^1_{s}-v^1_{s}\|^2_{{\cal L}{\cal B}^*} $.   The second term in the third line   is bounded by
\begin{equation*}
  C\| v_s^1\|^2_{\infty} \| u_s^1\|^2_1 \sup_{x\in \RR^d}\left| 
\hat{P}^1_{s,t}\varphi(x)-\check{P}^1_{s,t}\varphi(x)\right|^2
\leq  C
 \| v_s^1\|_{\infty}^2 \| u_s^1\|_1^2  \int_s^t \|u^1_r-v^1_r\|^2_{{\cal
L}{\cal B}^*}+ \|u^2_r-v^2_r\|^2_{{\cal
L}{\cal B}^*}dr
\end{equation*}
by  Lemma \ref{bornlip} b). The last threes lines are similarly dealt with and we  can easily conclude  as in Proposition \ref{uniqueness}.

\Box
 
 \section{Concluding remarks}
 
 We have developed models for dispersive and competitive multi species  population dynamics permitting nonlocal nonlinearity in the diffusive behavior of individuals. Depending on the value of the spatial competition range, the continuum (macro) limits of the individual (micro) dynamics turned out to be described by deterministic solutions of nonlocal cross-diffusion  systems with nonlocal or local  competition terms. These systems generalize usual diffusion-reaction systems  with nonlocal or local spatial nonlinearity in the diffusive coefficients, and  nonlocal or local nonlinearity in the reaction terms.
 These limiting objects can now be used as approximate objects for numerical simulation of spatial and ecological dynamics, when the individual behaviors depend  on the non homogeneous spatial densities of the different species. Of course, estimators of the relevant parameters of the phenomena under study have to be obtained first. 
 
 In the future it may also be worthwhile to elucidate  the situation where the spatial interaction range would be very small. This could thus justify by an individual-based approach the cross-diffusion models with local spatial interaction and local competition  extensively studied by the scientific community.

 \hspace{2cm}

 \textbf{Acknowledgements: }  The authors acknowledge  support of  ECOS-CONICYT C09E05 project. 
  J. Fontbona also thanks partial support from Basal-CONICYT grant  ``Center for Mathematical Modeling'' (CMM),   Millenium Nucleus  Stochastic Models of Disordered and Complex Systems NC120062 and the hospitality of \'Ecole Polytechnique.  S. M\'el\'eard  thanks  the Chair ``Mod\'elisation Math\'ematique et Biodiversit\'e '' of Veolia Environnement - \'Ecole Polytechnique - Museum National d'Histoire Naturelle - Fondation X,  and hospitality of CMM.  Both authors also thank Salom\'e Mart\'inez and Juan D\'avila for interesting and motivating discussions at the beginning of this work, as well as  for pointing out important references on cross-diffusion models.  We also thank the  anonymous referees for  pointing out some relevant references,  and for  remarks that allowed us to rectify or make more precise the statements of some of the hypothesis.

    \bibliography{biblio}

 \end{document}